\documentclass[10pt]{scrartcl}
\pdfoutput=1
\usepackage[labelfont={small},textfont={small}]{caption}
\usepackage{subcaption}
\usepackage{tocloft}
\usepackage{amssymb}
\usepackage{graphicx}
\usepackage{amsmath}
\usepackage[scaled=0.9]{helvet}
\usepackage{mathpazo}
\usepackage{paralist}
\usepackage{typearea}
\usepackage{upquote}
\usepackage[writefile]{listings}
\usepackage[dvipsnames]{xcolor}
\usepackage{microtype}
\usepackage{url}
\usepackage[pdftex,colorlinks]{hyperref}
\definecolor{darkblue}{cmyk}{1,1,0,0}
\definecolor{darkred}{cmyk}{0,1,0,0.7}
\hypersetup{anchorcolor=black,
  citecolor=darkblue, filecolor=darkblue,
  menucolor=darkblue,pagecolor=darkblue,urlcolor=darkblue,linkcolor=darkblue}

\title{{\DDEBIFCODE\ v. \version{}} Manual --- \\
  Bifurcation analysis
          of delay differential equations} 

\author{J. Sieber, K. Engelborghs, T. Luzyanina, G. Samaey, D. Roose}   

\date{\today}

\newcommand{\DDEBIFCODE}{\textsc{DDE-BIFTOOL}}
\newcommand{\ddebifweb}{\url{https://sourceforge.net/projects/ddebiftool}}
\newcommand{\ddebifarx}{\url{http://arxiv.org/abs/1406.7144}}
\newcommand{\ddebifwebold}{\url{http://twr.cs.kuleuven.be/research/software/delay/ddebiftool.shtml}}
\newcommand{\version}{3.1.1}

\newcommand{\file}[1]{\textbf{\texttt{#1}}}
\newcommand{\parm}[1]{\mathsf{#1}}
\newcommand{\define}[1]{\emph{#1}}
\newcommand{\demobase}{\url{../demos/index.html}}
\renewcommand{\i}{\mathrm{i}}
\newcommand{\T}{\mathrm{T}}
\renewcommand{\H}{\mathrm{H}}
\renewcommand{\d}{\mathrm{d}}
\newcommand{\RR}{\mathbb{R}}
\newcommand{\NN}{\mathbb{N}}
\newcommand{\ZZ}{\mathbb{Z}}
\newcommand{\CC}{\mathbb{C}}
\renewcommand{\Re}{\operatorname{Re}}
\renewcommand{\Im}{\operatorname{Im}}
\newcommand{\defeq}{:=}
\newcommand\numberthis{\addtocounter{equation}{1}\tag{\theequation}}

\newsavebox{\savepar}

\definecolor{var}{rgb}{0,0.28,0.28}
\definecolor{keyword}{rgb}{0,0,1}
\definecolor{comment}{rgb}{0,0.5,0}
\definecolor{string}{rgb}{0.6,0,0.5}
\definecolor{errmsg}{rgb}{1,0,0}
\newcommand{\blist}[1]{\mbox{\lstinline!#1!}}
\typearea{11}
\begin{document}
\pagenumbering{roman}
\maketitle



\noindent\textbf{\textsf{Keywords}} nonlinear dynamics, 
delay-differential equations, stability analysis, periodic solutions,
collocation methods, numerical bifurcation analysis, state-dependent
delay.






%
\renewcommand{\contentsname}{}
\tableofcontents
\clearpage
\pagenumbering{arabic}
\section{Citation, license, and obtaining the package}
\label{sec:app:get}
{\DDEBIFCODE} 
was started by Koen Engelborghs as part of his PhD at
the Computer Science Department of the K.U.Leuven under supervision
of Prof.\ Dirk Roose.

\paragraph{Citation}
Scientific publications for which the package DDE-BIFTOOL has been
used shall mention usage of the package DDE-BIFTOOL, and shall cite
the following publications to ensure proper attribution and
reproducibility:
\begin{compactitem}\item 
  K. Engelborghs, T. Luzyanina, and D. Roose. Numerical bifurcation
  analysis of delay differential equations using DDE-BIFTOOL, ACM
  Trans. Math. Softw. 28 (1), pp. 1-21, 2002.
\item\textbf{(this manual)} J. Sieber, K. Engelborghs, T. Luzyanina,
  G. Samaey, D. Roose . \DDEBIFCODE\ v. \version{} Manual ---
  Bifurcation analysis
          of delay differential equations, \ddebifarx{}.
\end{compactitem}
The implementation of normal forms for equilibria is based on
\begin{compactitem}
\item M.~M.~Bosschaert, B.~Wage and Y.~Kuznetsov: Description of the
  extension \texttt{ddebiftool\_nmfm}
  \url{http://ddebiftool.sourceforge.net/nmfm_extension_description.pdf}, 2015.
\item B. Wage: Normal form computations for Delay Differential
  Equations in DDE-BIFTOOL. Master Thesis, Utrecht University (NL),
  supervised by Y.A. Kuznetsov (\url{http://dspace.library.uu.nl/handle/1874/296912}, 2014.
\item M. M. Bosschaert: Switching from codimension 2 bifurcations of
  equilibria in delay differential equations. Master Thesis, Utrecht
  University (NL), supervised by Y.A. Kuznetsov,
  \url{http://dspace.library.uu.nl/handle/1874/334792}, 2016.
\end{compactitem}

All versions of this manual from v.~3.0 onward are available at
\ddebifarx{}.
\paragraph{License}
The following terms cover the use of the software package {\DDEBIFCODE}:
\bigskip\hrule
{\ttfamily
\begin{flushleft} {\parindent0pt
    \parskip5pt BSD 2-Clause license

Copyright (c) 2016, K.U. Leuven, Department of Computer Science,
K. Engelborghs, T. Luzyanina, G. Samaey. D. Roose, K. Verheyden,
J. Sieber, B. Wage, D. Pieroux

All rights reserved.

Redistribution and use in source and binary forms, with or without
modification, are permitted provided that the following conditions are
met:

1. Redistributions of source code must retain the above copyright
notice, this list of conditions and the following disclaimer.

2. Redistributions in binary form must reproduce the above copyright
notice, this list of conditions and the following disclaimer in the
documentation and/or other materials provided with the distribution.

THIS SOFTWARE IS PROVIDED BY THE COPYRIGHT HOLDERS AND CONTRIBUTORS
"AS IS" AND ANY EXPRESS OR IMPLIED WARRANTIES, INCLUDING, BUT NOT
LIMITED TO, THE IMPLIED WARRANTIES OF MERCHANTABILITY AND FITNESS FOR
A PARTICULAR PURPOSE ARE DISCLAIMED. IN NO EVENT SHALL THE COPYRIGHT
HOLDER OR CONTRIBUTORS BE LIABLE FOR ANY DIRECT, INDIRECT, INCIDENTAL,
SPECIAL, EXEMPLARY, OR CONSEQUENTIAL DAMAGES (INCLUDING, BUT NOT
LIMITED TO, PROCUREMENT OF SUBSTITUTE GOODS OR SERVICES; LOSS OF USE,
DATA, OR PROFITS; OR BUSINESS INTERRUPTION) HOWEVER CAUSED AND ON ANY
THEORY OF LIABILITY, WHETHER IN CONTRACT, STRICT LIABILITY, OR TORT
(INCLUDING NEGLIGENCE OR OTHERWISE) ARISING IN ANY WAY OUT OF THE USE
OF THIS SOFTWARE, EVEN IF ADVISED OF THE POSSIBILITY OF SUCH DAMAGE.
 }
\end{flushleft}
}
\hrule\bigskip

\paragraph{Download}

Upon acceptance of the above terms, one can obtain the package
{\DDEBIFCODE} (version \version{}) from
\begin{quote}
  \ddebifweb{}.
\end{quote}
Versions up to 3.0 and resources on theoretical background continue to
be available (under a different license) on
\begin{quote}
  \ddebifwebold{}.  
\end{quote}

\section{Version history}
\label{sec:changes}

\subsection{Changes from 3.1 to \version}
\label{sec:v31to311}
\begin{itemize}
\item Small bug fix: when changing focus on plot window during
  \blist{br_contn}, the online plotting no longer follows focus. This
  keeps the online bifurcation diagram in the same window. An optional
  named argument \blist{'plotaxis'} has been added to \blist{br_contn}
  to explicitly set the plot axes.
\item Added support for rotations (phase oscillators). Periodicity is
  enforced only up to multiples of $2\pi$. Demo
  \url{../demos/phase_oscillator/html/phase_oscillator.html} shows
  how one can track rotations. Demo was contributed by Azamat Yeldesbay.
\item Demos showing detection and computation of Bodganov-Takens
  bifurcation in
  \begin{quote}
    \url{../demos/Holling-Tanner/html/HollingTanner_demo.html}
  \end{quote}
  and cusp in
  \begin{quote}
    \url{../demos/cusp/html/cusp_demo.html}.
  \end{quote}
  Contributed by M.~M.~Boschaert and Y.~Kuznetsov.
\item A description of
  the mathematical formulas behind the normal form computations is now
  in \url{nmfm_extension_description.pdf}, by M.Bossschaert, B.~Wage,
  Y.~Kuznetsov.
\end{itemize}
\subsection{Changes from 3.0 to 3.1}
\label{sec:v3to31}
\paragraph{Change of License and move to Sourceforge}
D. Roose has permitted to change the license to a
Sourceforge-compliant BSD License. Thus, code and newest releases from
version 3.1 onward are now available from \ddebifweb{}. Older versions
will continue to be available from \ddebifwebold{}.
\paragraph{New feature: Normal form computation for bifurcations of
  equilibria}
The new functionality is only applicable for equations with constant
delay. Normal form coefficients can be computed through the extension
\texttt{ddebiftool\_extra\_nmfm}. This extension is included in the
standard \DDEBIFCODE{} archive, but the additional functions are kept
in a separate folder. The following bifurcations are currently
supported.
\begin{itemize}
\item Hopf bifurcation (coefficient $L_1$ determining criticality),
\item generalized Hopf (Bautin) bifurcation (of codimension two,
  typically encountered along Hopf curves)
\item Zero-Hopf interaction (Gavrilov-Guckenheimer bifurcation, of codimension two,
  typically encountered along Hopf curves)
\item Hopf-Hopf interaction (of codimension two, typically encountered
  along Hopf curves) 
\end{itemize}
The extension comes with a demo \texttt{nmfm\_demo}. The demos
\texttt{neuron}, \texttt{minimal\_demo} and \texttt{Mackey-Glass}
illustrate the new functionality, too. Background theory is given in \cite{W14}.
\subsection{Changes from 2.03 to 3.0}
\label{sec:v2to3}

\paragraph{New features}
\begin{itemize}
\item\label{int:pocont} \textbf{\textsf{Continuation of local periodic orbit
    bifurcations}} for systems with constant or state-dependent delay
  is now supported through the extension \texttt{ddebiftool\_extra\_psol}. This
  extension is included in the standard \DDEBIFCODE{} archive, but the
  additional functions are kept in a separate folder.
\item \textbf{\textsf{User-defined functions}} specifying the right-hand side
  and delays (such as \blist{sys_rhs} and \blist{sys_tau}) can have
  arbitrary names.  These user functions (with arbitrary names) get
  collected in a structure \blist{funcs}, which gets then passed on to
  the \DDEBIFCODE{} routines. This interface is similar to other
  functions acting on Matlab functions such as \blist{fzero} or
  \blist{ode45}. It enables users to add extensions such as
  \texttt{ddebiftool\_extra\_psol} without changing the core routines.
\item \textbf{\textsf{State-dependent delays}} can now have arbitrary levels of
  nesting (for example, periodic orbits of $\dot
  x(t)=\mu-x(t-x(t-x(t-x(t))))$ and their bifurcations can be
  tracked).
\item \textbf{\textsf{Vectorization}}\quad Continuation of periodic orbits
  and their bifurcations benefits (moderately) from
  vectorization of the user-defined functions.
\item \textbf{\textsf{Utilities}}\quad Recurring tasks (such as branching
  off at bifurcations, defining initial pieces of branches, or
  extracting the number of unstable eigenvalues) can now be performed
  more conveniently with some auxiliary functions provided in a
  separate folder \texttt{ddebiftool\_utilities}. See
  \url{../demos/neuron/html/demo1_simple.html} for a demonstration.
\item \textbf{\textsf{Bugs fixed}}\quad Some bugs and problems have been
  fixed in the implementation of the heuristics applied to choose the
  stepsize in the computation of eigenvalues of equilibria
  \cite{VLR08}.
\item \textbf{\textsf{Continuation of relative equilibiria and
      relative periodic orbits}} and their local bifurcations for
  systems with constant delay and rotational symmetry (saddle-node
  bifurcation, Hopf bifurcation, period-doubling, and torus
  bifurcation) is now supported through the extension
  \texttt{ddebiftool\_extra\_rotsym}.
\end{itemize}
Extensions come with demos and separate documentation.

\paragraph{Change of user interface}
Versions from 3.0 onward have a different the user interface for many
\DDEBIFCODE{} functions than versions up to 2.03. They add one
additional input argument \blist{funcs} (this new argument comes
always \emph{first}).  Since \DDEBIFCODE{} v.\,3.0 changed the user
interface, scripts written for \DDEBIFCODE{} v.\,2.0x or earlier will
not work with versions later than 3.0. For this reason version 2.03
will continue to be available. Users of both versions should ensure
that only one version is in the Matlab path at any time to avoid
naming conflicts.

\section{Capabilities and related reading and
  software}
\label{sec:intro}
 {\DDEBIFCODE} consists of a set of routines running in
Matlab\footnote{\url{http://www.mathworks.com}} \cite{Mat00} or
octave\footnote{\url{http://www.gnu.org/software/octave}}, both widely
used environments for scientific computing.  The aim of the package is
to provide a tool for numerical bifurcation analysis of steady state
solutions and periodic solutions of differential equations with
constant delays (called DDEs) or state-dependent delays (here called
sd-DDEs).  It also allows users to compute homoclinic and heteroclinic
orbits in DDEs (with constant delays).  

\paragraph{Capabilities}
{\DDEBIFCODE} can perform the following computations:
\begin{compactitem}
\item continuation of steady state solutions (typically in a single
  parameter);
\item approximation of the rightmost, stability-determining roots of
  the characteristic equation which can further be corrected using a
  Newton iteration;
\item continuation of steady state folds and Hopf bifurcations
  (typically in two system parameters);
\item continuation of periodic orbits using orthogonal collocation
  with adaptive mesh selection (starting from a previously computed
  Hopf point or an initial guess of a periodic solution profile);
\item approximation of the largest stability-determining Floquet
  multipliers of periodic orbits;
\item branching onto the secondary branch of periodic solutions at a
  period doubling bifurcation or a branch point;
\item continuation of folds, period doublings and torus bifurcations
  (typically in two system parameters) using the extension
  \texttt{ddebiftool\_extra\_psol};
\item computation of normal form coefficients for Hopf bifurcations
  and codimension-two bifurcations along Hopf bifurcation curves
  (typically in two system parameters) using the extension
  \texttt{ddebiftool\_extra\_nmfm};
\item continuation of connecting orbits (using the appropriate number
  of parameters).  
\end{compactitem}
All computations can be performed for problems with an arbitrary
number of discrete delays. These delays can be either parameters or
functions of the state.
(The only exception are computations of connecting orbits, which
support only problems with delays as parameters at the moment.)

A practical difference to AUTO or MatCont is that the package does not detect
bifurcations automatically because the computation of eigenvalues or
Floquet multipliers may require more computational effort than the
computation of the equilibria or periodic orbits (for example, if the
system dimension is small but one delay is large). Instead the
evolution of the eigenvalues can be computed along solution branches
in a separate step if required. This allows the user to detect and
identify bifurcations.



\paragraph{About this manual --- related reading}
This manual documents version~\version. Earlier versions of the manual
for earlier versions of \DDEBIFCODE{} continue to be available at the
web addresses
\begin{tabbing}
  versions $\geq3.0$\qquad \= \ddebifarx{}\\
  versions $\leq2.03$ \> \ddebifwebold{}.
\end{tabbing}
For readers who intend to analyse only systems with constant delays,
the parts of the manual related to systems with state-dependent delays
can be skipped (sections \ref{sd_dde}, \ref{sys_def2}).  In the rest
of this manual we assume the reader is familiar with the notion of a
delay differential equation and with the basic concepts of bifurcation
analysis for ordinary differential equations.  The theory on delay
differential equations and a large number of examples are described in
several books. Most notably the early
\cite{Bell63,Driv77,El's73,Hale77a,Kolm86} and the more recent
\cite{Azbe91,Kolm92,Hale93,Diek95,Kolm99}.  Several excellent books
contain introductions to dynamical systems and bifurcation theory of
ordinary differential equations, see, e.g.,
\cite{Argy94,Chow82,Guck83,Kuzn04,Seyd94}.

\paragraph{Turorial demos}
The tutorial demos \texttt{demo1} and
\texttt{sd\_demo}, providing a step-by-step walk-through for the
typical working mode with \DDEBIFCODE{} are included as separate
\texttt{html} files, published directly from the comments in the demo
code. See \demobase{} for links to all demos, many of which are
extensively commented.

\paragraph{Related software}
A large number of packages exist for numerical continuation and
bifurcation analysis of systems of ordinary differential
equations. Currently maintained packages are
\begin{tabbing}
  \texttt{AUTO}\qquad
  \= url: \url{http://sourceforge.net/projects/auto-07p}\qquad\=
  using FORTRAN or C \cite{Doed99,Doed07},\\
  \texttt{MatCont}\>url: \url{http://sourceforge.net/projects/matcont/}\>
  for Matlab \cite{DGK03,G00}, and\\
  \texttt{Coco}\> url: \url{http://sourceforge.net/projects/cocotools}\> for
  Matlab \cite{DS13}.\\[1ex]
  For delay differential equations the package\\[0.5ex]
  \texttt{knut}\> url: \url{http://gitorious.org/knut}\> using C++
\end{tabbing}
(formerly \texttt{PDDECONT}) is available as a
stand-alone package (written in C++, but with a user interface
requiring no programming). This package was developed in parallel with
\DDEBIFCODE{} but independently by R.~Szalai \cite{SSH06,RS07}.  For
simulation (time integration) of delay differential equations the
reader is, e.g., referred to the packages ARCHI, DKLAG6, XPPAUT,
DDVERK, RADAR and dde23, see
\cite{Paul95,Thom97,Erme98,Enri97,Sham00,Gugl07}.  Of these, only
XPPAUT has a graphical interface (and allows limited stability
analysis of steady state solutions of DDEs along the lines of
\cite{Luzy96}). TRACE-DDE is a Matlab tool (with graphical interface)
for linear stability analysis of linear constant-coefficient DDEs
\cite{breda09}.


\begin{figure}[t]
\begin{center}
\includegraphics[width=0.8\textwidth]{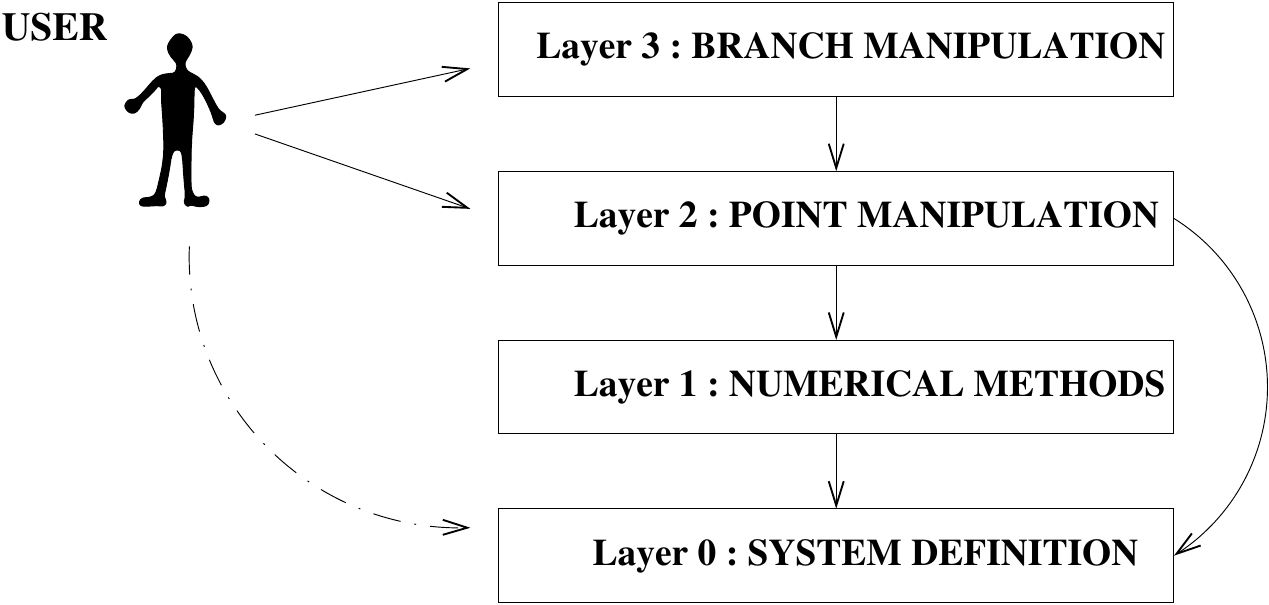}
\end{center}
\caption{\label{struct_pic}
The structure of {\DDEBIFCODE}. Arrows indicate the 
calling ($-$) or writing ($\cdot-$) of routines in a certain layer.} 
\end{figure}

\section{Structure of {\DDEBIFCODE}}\label{code_struct}

The structure of the package is depicted in figure \ref{struct_pic}.
It consists of four layers. 

Layer 0 contains the system definition and consists of routines which
allow to evaluate the right hand side $f$ and its derivatives,
state-dependent delays and their derivatives and to set or get the
parameters and the constant delays.  It should be provided by the user
and is explained in more detail in section \ref{sec:system:def}.  All
user-provided functions are collected in a single structure (called
\blist{funcs} in this manual), and are passed on by the user as
arguments to layer-3 or layer-2 functions. \textbf{\emph{Note that
    this is a change in user interface between version 2.03 and
    version~3.0!}}

Layer 1 forms the numerical core of the package and is (normally)
not directly accessed by the user. The numerical methods used
are explained
briefly in section \ref{numerical_methods}, more
details can be found in the papers 
\cite{Luzy96,Enge99a,Enge99b,en_d01,engel01,luz01,homoclinic}
and in \cite{Enge00}. Its functionality is hidden by and used
through layers 2 and 3.

Layer 2 contains routines to manipulate individual points.  Names of
routines in this layer start with "\file{p\_}".  A point has one of
the following five types.  It can be a steady state point (abbreviated
\blist{'stst'}), steady state Hopf (abbreviated \blist{'hopf'}) or
fold (abbreviated \blist{'fold'}) bifurcation point, a periodic solution point
(abbreviated \blist{'psol'}) or a connecting orbit point (abbreviated
\blist{'hcli'}). Furthermore, a point can contain additional information
concerning its stability.  Routines are provided to compute individual
points, to compute and plot their stability and to convert points from
one type to another.

Layer 3 contains routines to manipulate branches.  Names of routines
in this layer start with "\file{br\_}". A branch is structure
containing an array of (at least two) points, three sets of method
parameters and specifications concerning the free parameters.  The
\blist{'point'} field of a branch contains an array of points of the
same type ordered along the branch.  The \blist{'method'} field
contains parameters of the computation of individual points, the
continuation strategy and the computation of stability.  The
\blist{'parameter'} field contains specification
of the free parameters (which are allowed to vary along the branch),
parameter bounds and maximal step sizes.  Routines are provided to
extend a given branch (that is, to compute extra points using
continuation), to (re)compute stability along the branch and to
visualize the branch and/or its stability.

Layers 2 and 3 require specific data structures, explained in
section~\ref{data_structures}, to represent points, stability
information, branches, to pass method parameters and to specify
plotting information.  Usage of these layers is
demonstrated 
through a step-by-step analysis of the demo systems \texttt{neuron},
\texttt{sd\_demo} and \texttt{hom\_demo} (see
\demobase{}).  Descriptions of input/output parameters and
functionality of all routines in layers 2 and 3 are given in
sections~\ref{point_manipulation} respectively
\ref{branch_manipulation}.

\section{Delay differential equations}\label{explain_dde}
\label{sec:ddes}
This section introduces the mathematical notation that we refer to in
this manual to describe the problems solved by \DDEBIFCODE{}.
\subsection{Equations with constant delays}\label{dde}

Consider the system of delay differential equations with constant
delays (DDEs),
\begin{equation}\label{the_dde_type}
\frac{\d}{\d t}{x(t)}=f(x(t),x(t-\tau_1),\ldots,x(t-\tau_m),\eta),
\end{equation}
where $x(t)\in\RR^n$, $f:\RR^{n(m+1)}\times\RR^p
\rightarrow\RR^n$ is a nonlinear smooth function
depending on a number of parameters $\eta\in\RR^p$, and delays
$\tau_i>0$, $i=1,\ldots,m$.
Call $\tau$ the maximal delay,
\[
\tau=\max_{i=1,\ldots,m}\tau_i.
\]
The linearization of (\ref{the_dde_type}) around a solution $x^*(t)$ 
is the \define{variational equation}, given by,
\begin{equation}\label{the_var_equa}
\frac{\d}{\d t}{y(t)}=A_0(t)y(t)+\sum_{i=1}^m A_i(t)y(t-\tau_i),
\end{equation}
where, using $f\equiv f(x^0,x^1,\ldots,x^m,\eta)$,  
\begin{equation}\label{A_def}
A_i(t)=\frac{\partial f}{\partial x^i}(x^*(t),x^*(t-\tau_1),\ldots,x^*(t-\tau_m),\eta), 
\ i=0,\ldots,m. 
\end{equation}

\subsubsection{Steady states}
\label{sec:dde:stst}
If $x^*(t)$ corresponds to a steady state solution,
\[
x^*(t)\equiv x^*\in\RR^n,\mathrm{\ with\ }f(x^*,x^*,\ldots,x^*,\eta)=0,
\]
then the matrices 
$A_i(t)$ are constant, $A_i(t)\equiv A_i$, and the corresponding 
variational equation (\ref{the_var_equa})
leads to a \define{characteristic equation}. Define the $n\times n$-dimensional
matrix $\Delta$ as
\begin{equation}
  \Delta(\lambda)=\lambda I - A_0 - \sum_{i=1}^m A_i e^{-\lambda\tau_i}.
\label{eq:deltadef}
\end{equation}
Then the characteristic equation reads,
\begin{equation}\label{the_char_eq}
\det(\Delta(\lambda))=0.
\end{equation}
Equation (\ref{the_char_eq}) has an infinite number of 
roots $\lambda\in\CC$ which determine the stability of the steady state
solution $x^*$.
The steady state solution is (asymptotically) stable provided all
roots of the characteristic equation (\ref{the_char_eq}) have
negative real part; it is unstable if there exists a root with positive
real part.
It is known that the number of roots in any right half plane
$\Re(\lambda)>\gamma$, $\gamma\in\RR$ is finite, hence, the
stability is always determined by a finite number of roots.

Bifurcations occur whenever roots move through the imaginary
axis as one or more parameters are changed.
Generically a fold bifurcation (or turning point) occurs when
the root is real (that is, equal to zero) and a 
Hopf bifurcation occurs when a pair of complex conjugate roots crosses the imaginary axis.

\subsubsection{Periodic orbits}
\label{sec:dde:psol}
A periodic solution $x^*(t)$ is a solution which repeats itself after
a finite time, that is,
\[ 
x^*(t+T)=x^*(t),\mathrm{\ for\ all\ }t. 
\]
Here $T>0$ is the period.  The stability around the periodic solution
is determined by the time integration operator $S(T,0)$ which
integrates the variational equation (\ref{the_var_equa}) around
$x^*(t)$ from time $t=0$ over the period.  This operator is called the
\define{monodromy operator} and its (infinite number of) eigenvalues,
which are independent of the starting moment $t=0$, are called the
\define{Floquet multipliers}.  Furthermore, if $S(T,0)^k$ is compact
for $k>\tau/T$. Thus, there are at most finitely many Floquet
multipliers outside of any ball around the origin of the complex
plane.

For autonomous systems there is always a \define{trivial} Floquet
multiplier at unity, corresponding to a perturbation along the time
derivative of the periodic solution. The periodic solution is
exponentially stable provided all multipliers (except the trivial one)
have modulus smaller than unity, it is exponenially unstable if there
exists a multiplier with modulus larger than unity.

\subsubsection{Connecting orbits}
\label{sec:dde:hcli}
We call a solution $x^*(t)$ of (\ref{the_dde_type}) at $\eta=\eta^*$ a 
\textit{connecting orbit} if the limits 
\begin{equation}
\lim_{t\to -\infty} x^*(t)=x^{-}, \qquad \lim_{t\to +\infty} x^*(t)=x^{+},
\end{equation}
exist.  For continuous $f$, $x^-$ and $x^+$ 
are steady state solutions. 
If $x^-=x^+$, the orbit is called homoclinic, otherwise it is heteroclinic. 

\subsection{Equations with state-dependent delays}\label{sd_dde}

Consider the system of delay differential equations with
state-dependent delays (sd-DDEs),
\begin{equation}\label{the_dde_type2}
\left\{
\begin{aligned}
\frac{\d}{\d t}{x(t)}&=f(x_0,x_1,\ldots,x_m,\eta),\\
x_j&=x(t-\tau_j(x_0,\ldots,x_{j-1},\eta)\quad\mbox{($\tau_0=0$, $j=1,\ldots,m$)}\mbox{,}
\end{aligned}
\right.
\end{equation}
where $x(t)\in\RR^n$, and
\begin{align*}
 f:&\RR^{n(m+1)}\times\RR^p \to\RR^n\\
 \tau_j:&\RR^{n\,j}\times\RR^p\to[0,\infty)
\end{align*}
are smooth functions depending of their arguments. The right-hand side
$f$ depends on $m+1$ states $x_j=x(t-\tau_j)\in\RR^n$ ($j=0,\ldots,m$)
and $p$ parameters $\eta\in\RR^p$. The $j$th delay function $\tau_j$
depends on all previously defined $j-1$ states $x(t-\tau_i)\in\RR^n$
($i=0,\ldots,j-1$) and $p$ parameters $\eta\in\RR^p$. This
definition permits the user to formulate sd-DDEs with arbitrary levels
of nesting in their function arguments.

The linearization around a solution $(x^*(t),\eta^*)$ of
\eqref{the_dde_type2} (the \define{variational equation}) with respect
to $x$ is given by (see \cite{HKWW06}, we are using the notation
$x^*_0=x^*(t)$, $\tau_j^*(t)=\tau_j(x^*_0,\ldots,x^*_{j-1},\eta^*)$ and
$x^*_j=x^*(t-\tau^*_j(t))$ for $j\geq1$)
\begin{equation}\label{the_var_equa2}
\begin{aligned}
\frac{\d}{\d t}{y(t)}&=\sum_{j=0}^mA_j(t)Y_k \\
Y_0&=y(t)\\
Y_j&=y(t-\tau^*_j(t))-(x^*)'(t-\tau_j(t))\sum_{k=0}^{j-1}B_{k,j}(t)Y_k\mbox{,\quad($j=1\ldots m$)}
\end{aligned}
\end{equation}
where $(x^*)^{'}(t)={\d}x^*(t)/{\d}t$, and
\begin{equation}\label{A_def2}
\begin{aligned}
  A_j(t)&=\frac{\partial f}{\partial x^j}(x^*_0,x^*_1,\ldots,x^*_m,\eta)\in\RR^{n\times n}
  \mbox{,\quad ($j=0,\ldots,m$),}\\
  B_{j,k}(t)&=\frac{\partial \tau_j}{\partial x^k}(x^*_0,x^*_1,\ldots,x^*_{j-1},\eta^*)\in\RR^{1\times n}\mbox{,\quad ($k=0,\ldots,j-1$, $j=1,\ldots,m$).} 
\end{aligned}
\end{equation}

If $(x^*(t),\tilde{\tau}^*(t))$ corresponds to a steady state
solution, then $x^*(t)=x^*_0=\ldots=x^*_m\equiv x^*\in\RR^n$, and
$\tau^*_j(t)\equiv \tau_j(x^*,\ldots,x^*,\eta^*)$ for all $j\geq1$, with
\[
f(x^*,x^*,\ldots,x^*,\eta^*)=0
\]
then the matrices $A_i(t)$ are constant, $A_i(t)\equiv A_i$, and the
vectors $B_{i,j}(t)$ consist of zero elements only.  In this case, the
corresponding variational equation \eqref{the_var_equa2} is a constant
delay differential equation and it leads to the characteristic
equation (\ref{the_char_eq}), i.e.~a characteristic equation with
constant delays. Hence the stability analysis of a
steady state solution of \eqref{the_dde_type2} is similar to the
stability analysis of \eqref{the_dde_type}.

Note that the right-hand side $f$, when considered as a functional
mapping a history segment into $\RR^n$ is not locally Lipschitz
continuous. This creates technical difficulties when considering an
sd-DDE of type \eqref{the_dde_type2} as an infinite-dimensional
system, because the solution does not depend smoothly on the
initial condition (see Hartung \emph{et al }\cite{HKWW06} for a
detailed review). However, periodic boundary-value problems for
\eqref{the_dde_type2} can be reduced to finite-dimensional systems of
algebraic equations that are as smooth as the coefficient functions
$f$ and $\tau_j$ \cite{S12}. This implies that all periodic orbits and
their bifurcations and stability as computed by \DDEBIFCODE{} behave as
expected. In particular, branching off at Hopf bifurcations and period
doubling works in the same way as for constant delays (the proof for the Hopf
bifurcation in sd-DDEs is also given in \cite{S12}).


Moreover, Mallet-Paret and Nussbaum \cite{MN11} proved that the
stability of the linear variational equation \eqref{the_var_equa2}
indeed reflects the {\it local stability} of the solution
$(x^*(t),\tilde{\tau}^*(t))$ of \eqref{the_dde_type2}. 
For details on the relevant theory and numerical bifurcation analysis
of differential equations with state-dependent delay see
\cite{luz01,HKWW06} and the references therein.

\section{System definition}\label{sec:system:def}
\subsection{Change from \DDEBIFCODE{} v.~2.03 to
  v.~3.0} Note that in \DDEBIFCODE{} \version{} all
user-provided functions can be arbitrary function handles, collected
into a structure using the function \blist{set_funcs}, explained in
section~\ref{sec:funcs}. The only typical mandatory functions for the
user to provide are the right-hand side (\blist{'sys_rhs'}) and the
function returning the delay indices (\blist{'sys_tau'}). The names of
the user functions can be arbitrary, and user functions can be
anonymous. This is a change from previous versions. See the tutorials
in \demobase{} for examples of usage, and the function description in
section~\ref{sec:funcs} for details.
\subsection{Equations with constant delays}\label{sys_def1}

As an illustrative example we will use the following system of delay 
differential
equations, taken from \cite{Shay99},
\begin{equation}\label{example_sys}
\left\{
\begin{array}{l}
\dot{x_1}(t)=-\kappa x_1(t)+\beta \tanh(x_1(t-\tau_s))+a_{12}\tanh(x_2(t-\tau_2)) \\
\dot{x_2}(t)=-\kappa x_2(t)+\beta \tanh(x_2(t-\tau_s))+a_{21}\tanh(x_1(t-\tau_1)) .
\end{array}
\right.
\end{equation}
This system models two coupled neurons with time delayed connections.
It has two components ($x_1$ and $x_2$), three delays ($\tau_1$,
$\tau_2$ and $\tau_s$), and four parameters ($\kappa$, $\beta$,
$a_{12}$ and $a_{21}$).  The demo \texttt{neuron} (see
\demobase{}) walks through the
bifurcation analysis of system \eqref{example_sys} step by step to
demonstrate the working pattern for \DDEBIFCODE.

To define a system, the user should provide the following Matlab
functions, given in the following paragraphs for system \eqref{example_sys}.

\subsubsection{Right-hand side --- \blist{sys_rhs}}\label{sec:constrhs} 
The right-hand side is a function of two arguments. For our example
\eqref{example_sys}, this would have the form (giving the right-hand
side the name \blist{neuron_sys_rhs})
\begin{lstlisting}[frame=lines,label=neuron_sys_rhs,caption={Definition for right-hand side of \eqref{example_sys} as a variable.}]
neuron_sys_rhs=@(xx,par)[...
  -par(1)*xx(1,1)+par(2)*tanh(xx(1,4))+par(3)*tanh(xx(2,3));...
  -par(1)*xx(2,1)+par(2)*tanh(xx(2,4))+par(4)*tanh(xx(1,2))];  
%par=[\kappa,\beta, a_{12}, a_{21},\tau_1,\tau_2, \tau_s]
\end{lstlisting}
Meaning of the arguments of the right-hand side function:
\begin{itemize}
\item $\blist{xx}\in\RR^{n\times (m+1)}$ contains the state
variable(s) at the present and in the past,
\item $\blist{par}\in\RR^{1\times
  p}$ contains the parameters, $\blist{par}=\eta$.
\end{itemize}
The delays $\tau_i$ ($i=1\ldots,m$) are considered to be part of the
parameters ($\tau_i=\eta_{j(i)}$, $i=1,\ldots,m$).  This is natural
since the stability of steady solutions and the position and stability
of periodic solutions depend on the values of the delays.  Furthermore
delays can occur both as a `physical' parameter and as delay, as in
$\dot{x}=\tau x(t-\tau)$.  From these inputs the right hand side $f$
is evaluated at time $t$. Notice that the parameters have a specific
order in $\parm{par}$ indicated in the comment line.

An alternative (vectorized) form would be
\begin{lstlisting}[frame=lines,label=neuron_sys_rhs_vec,caption={Alternative definition of the right-hand side of \eqref{example_sys}, vectorized for speed-up of periodic orbit computations.}]
  neuron_sys_rhs=@(xx,p)[...
  -p(1)*xx(1,1,:)+p(2)*tanh(xx(1,4,:))+p(3)*tanh(xx(2,3,:));....
  -p(1)*xx(2,1,:)+p(2)*tanh(xx(2,4,:))+p(4)*tanh(xx(1,2,:))];
\end{lstlisting}
Note the additional colon in argument \blist{xx} and compare to
Listing~\ref{neuron_sys_rhs}. The form shown in
Listing~\ref{neuron_sys_rhs_vec} can be called in many points along a
mesh simultaneously, speeding up the computations during analysis of
periodic orbits.

\subsubsection{Delays --- \blist{sys_tau}} \label{sec:consttau}
For constant delays another function is required which returns the
\emph{position} of the delays in the parameter list. For our example,
this is
\begin{lstlisting}
  neuron_tau=@()[5 6 7];
\end{lstlisting}
This function has no arguments for constant delays, and returns a row
vector of indices into the parameter vector.

\subsubsection{Jacobians of right-hand side --- \blist{sys_deri} (optional, but recommended)}\label{sec:constjac}
Several derivatives of the right hand side function $f$ need to be
evaluated during bifurcation analysis. By default, \DDEBIFCODE{} uses
a finite-difference approximation, implemented in
\file{df\_deriv.m}. For speed-up or in case of convergence
difficulties the user may provide the Jacobians of the right-hand side
analytically as a separate function. Its header is of the format
\begin{lstlisting}
  function J=sys_deri(xx,par,nx,np,v)
\end{lstlisting}
Arguments:
\begin{itemize}
\item $\blist{xx}\in\RR^{n\times (m+1)}$ contains the state
variable(s) at the present and in the past (as for the right-hand side);
\item $\blist{par}\in\RR^{1\times
  p}$ contains the parameters, $\blist{par}=\eta$  (as for the right-hand side);
\item \blist{nx} (empty, one integer or two integers) index (indices) of
  \blist{xx} with respect to which the right-hand side is to be
  differentiated
\item \blist{np} (empty or integer) whether right-hand side is to be
  differentiated with respect to parameters
\item \blist{v} (empty or $\CC^n$) for mixed derivatives with respect
  to \blist{xx}, only the product of the mixed derivative with
  \blist{v} is needed.
\end{itemize}
The result \blist{J} is a matrix of
partial derivatives of $f$ which depends on the type of derivative
requested via \blist{nx} and \blist{np} multiplied with \blist{v} (when
nonempty), see table \ref{deri_requested}.
\begin{table}[htbp]
\begin{center}
\begin{tabular}{ccc@{\hspace*{5em}}l}
\noalign{\medskip}\hline\noalign{\smallskip}\blist{length(nx)} & \blist{length(np)}  & \blist{v} & \blist{J} 
\\\noalign{\smallskip}\hline\noalign{\medskip}
1         & 0         & empty      & 
$\frac{\textstyle\partial f}{\textstyle\partial x^{\blist{nx(1)}}}
=A_{\blist{nx(1)}}\in\RR^{n\times n}$ \\[3ex]
0         & 1         & empty      & 
$\frac{\textstyle\partial f}{\textstyle\partial \eta_{\blist{np(1)}}}
\in\RR^{n\times 1}$ \\[3ex]
1         & 1         & empty      & 
$\frac{\textstyle\partial^2 f}{\textstyle\partial x^{\blist{nx(1)}}\partial \eta_{\blist{np(1)}}}
\in\RR^{n\times n}$ \\
2         & 0         & $\in\CC^{n\times1}$ & 
$\frac{\textstyle\partial}{\textstyle\partial x^{\blist{nx(2)}}}
\left(A_{\blist{nx(1)}}v\right)
\in\CC^{n\times n}$\\\noalign{\medskip}\hline
\end{tabular}
\caption{\label{deri_requested}
Results of the function \file{sys\_deri} depending on its
input parameters \blist{nx}, \blist{np} and \blist{v}
using $f\equiv f(x^0,x^1,\ldots,x^m,\eta)$.}
\end{center}
\end{table}

\blist{J} is defined as follows. Initialize \blist{J} with $f$. If
\blist{nx} is nonempty take the derivative of \blist{J} with respect
to those arguments listed in \blist{nx}'s entries. Each entry of
\blist{nx} is a number between $0$ and $m$ based on $f\equiv
f(x^0,x^1,\ldots,x^m,\eta)$.  E.g., if \blist{nx} has only one element
take the derivative with respect to $x^{\blist{nx(1)}}$.  If it has
two elements, take, of the result, the derivative with respect to
$x^{\blist{nx(2)}}$ and so on.  Similarly, if \blist{np} is nonempty
take, of the resulting \blist{J}, the derivative with respect to
$\eta_{\blist{np(i)}}$ where $i$ ranges over all the elements of
\blist{np}, $1\leq i \leq p$.  Finally, if $v$ is not an empty vector
multiply the result with $v$.  The latter is used to prevent \blist{J}
from being a tensor if two derivatives with respect to state variables
are taken (when \blist{nx} contains two elements).  Not all possible
combinations of these derivatives have to be provided.  In the current
version, \blist{nx} has at most two elements and \blist{np} at most
one.  The possibilities are further restricted as listed in table
\ref{deri_requested}.

In the last row of table \ref{deri_requested} the elements of \blist{J}
are given by,
\[
\blist{J}_{i,j}=\left[\frac{\partial}{\partial x^{\blist{nx(2)}}}
A_{\blist{nx(1)}}v\right]_{i,j}
=\frac{\partial}{\partial x_j^{\blist{nx(2)}}}
\left(\sum_{k=1}^n\frac{\partial f_i}{\partial x_k^{\parm{nx(1)}}} v_k
\right),
\]
with $A_l$ as defined in \eqref{A_def}.

The resulting routine is quite long, even for the small system
\eqref{example_sys}; see Listing~\ref{neuron:sys:deri} in
Appendix~\ref{sec:sys:deri} for a printout of the function body.  Furthermore,
implementing so many derivatives is an activity prone to a number of
typing mistakes. Hence a default routine \blist{df\_deriv} is
available which implements finite difference formulas to approximate
the requested derivatives (using several calls to the right-hand
side. It is, however, recommended to provide at least the first order
derivatives with respect to the state variables using analytical
formulas. These derivatives occur in the determining systems for fold
and Hopf bifurcations
and for connecting orbits,  
and in the computation of characteristic roots and Floquet multipliers.
All other derivatives are only necessary in the 
Jacobians of the respective Newton procedures and thus
influence only the convergence speed.

\subsection{Equations with state-dependent
  delays}\label{sys_def2}
\DDEBIFCODE{} also permits the delays to depend on parameters and the
state. If at least one delay is state-dependent then the format and
semantics of the function specifying the delays, \blist{sys_tau}, is
different from the format used for constant delays in
section~\ref{sec:consttau} (it now provides the \emph{values} of the
delays).  Note that for a system with only constant delays we
recommend the use of the system definitions as described in
section~\ref{sys_def1} to reduce the computational effort.

As an illustrative example we will use the following system of delay 
differential equations,
\begin{equation}\label{example_sys2}
  \begin{split}
    \frac{\d}{\d t}x_1(t)&=\frac{1}{p_1+x_2(t)}\left(1-p_2x_1(t)x_1(t-\tau_3)
      x_3(t-\tau_3)+p_3x_1(t-\tau_1)x_2(t-\tau_2)\right),\\
    \frac{\d}{\d t}x_2(t)&=\frac{p_4 x_1(t)}{p_1+x_2(t)}+
    p_5\tanh(x_2(t-\tau_5))-1,\\
    \frac{\d}{\d t}x_3(t)&=p_6(x_2(t)-x_3(t))-p_7(x_1(t-\tau_6)-x_2(t-\tau_4))e^{-p_8 \tau_5},\\
    \frac{\d}{\d t}x_4(t)&=x_1(t-\tau_4)e^{-p_1 \tau_5} -0.1,\\
    \frac{\d}{\d t}x_5(t)&=3(x_1(t-\tau_2)-x_5(t))-p_9,
  \end{split}
\end{equation}
where
\begin{align*}
\tau_1, \tau_2 &\mbox{ are constant delays},\\
\tau_3&=2+p_5\tau_1x_2(t)x_2(t-\tau_1),\\
\tau_4&=1-\frac{1}{1+x_1(t)x_2(t-\tau_2)},\\
\tau_5&=x_4(t),\\
\tau_6&=x_5(t).
\end{align*}

This system has five components $(x_1,\ldots,x_5)$, six delays
$(\tau_1,\ldots,\tau_6)$ and eleven parameters $(p_1,\ldots,p_{11})$,
where $p_{10}=\tau_1$ and $p_{11}=\tau_2$.  A step-by-step tutorial
for analysis of sd-DDEs is given in demo \texttt{sd\_demo} (see
\demobase{}) of this system using \eqref{example_sys2}.

To define a system with state-dependent delays, the user should
provide the following Matlab functions, given in the following
sections for system \eqref{example_sys2}.

\subsubsection{Right-hand side --- \blist{sys_rhs}}
\label{sec:sdrhs}
The definition and functionality of this routine is equivalent to the
one described in section \ref{sec:constrhs}.  Notice that the argument
\blist{xx} contains the state variable(s) at the present and in the
past, $\blist{xx}=[x(t)\ x(t-\tau_1)\ \ldots\ \ldots\ x(t-\tau_m)]$.
Possible constant delays ($\tau_1$ and $\tau_2$ in example
\eqref{example_sys2}) are also considered to be part of the
parameters. See Listing \ref{sd_rhs} for the right-hand side to be
provided for example \eqref{example_sys2}.
\begin{lstlisting}[frame=lines,float,label=sd_rhs,
basicstyle={\ttfamily\small},caption={Listing of right-hand side
  \blist{'sys_rhs'} function for \eqref{example_sys2} function for
  \eqref{example_sys2}, here called \blist{sd_rhs}.}]
function f=sd_rhs(xx,par)
f(1,1)=(1/(par(1)+xx(2,1)))*(1-par(2)*xx(1,1)*xx(1,4)*xx(3,4)+...
        par(3)*xx(1,2)*xx(2,3));
f(2,1)=par(4)*xx(1,1)/(par(1)+xx(2,1))+par(5)*tanh(xx(2,6))-1;
f(3,1)=par(6)*(xx(2,1)-xx(3,1))-par(7)*(xx(1,7)-...
        xx(2,5))*exp(-par(8)*xx(4,1));
f(4,1)=xx(1,5)*exp(-par(1)*xx(4,1))-0.1;
f(5,1)=3*(xx(1,3)-xx(5,1))-par(9); 
end



\end{lstlisting}

\subsubsection{Delays --- \blist{sys_tau} and \blist{sys_ntau}}
\label{sec:sdtau}
The format and semantics of the routines specifying the delays differ
from the one described in section~\ref{sec:consttau}. The user has to provide two functions:
\begin{lstlisting}
  function ntau=sys_ntau()        
  function  tau=sys_tau (ind,xx,par)
\end{lstlisting}
The function \blist{sys_ntau} has no arguments and  returns the number of (constant and state-dependent) delays. For the example \eqref{example_sys2}, this could be the anonymous function \blist{@()6;}.

The function \blist{sys_tau} has the three arguments:
\begin{itemize}
\item \blist{ind} (integer $\geq1$) indicates, which delay is to be returned;
\item \blist{xx} ($ n\times\blist{ind}$-matrix) is the state:
  \blist{xx(:,1)}$=x(t)$, \blist{xx(:,k)}$=x(t-\tau_{k-1})$ for
  $k=2\ldots\blist{ind}$;
\item \blist{par} (row vector) is the vector of system parameters
\end{itemize}
The output is the value of $\tau_{\blist{ind}}$, the \blist{ind}'th
delay. \DDEBIFCODE{} calls \blist{sys_tau} $m$ times
where $m=$\blist{sys_ntau()}. In the first call \blist{xx} is the
$n\times1$ vector, $x(t)$ such that $\tau_1$ may depend on $x(t)$ and
\blist{par}. After the first call to \blist{sys_tau}, \DDEBIFCODE{}
computes $x(t-\tau_1)$. In the second call to \blist{sys_tau}
\blist{xx} is a $n\times2$ matrix, consisting of $[x(t),x(t-\tau_1)]$
such that the delay $\tau_2$ may depend on $x(t)$ , $x(t-\tau_1)$ and
\blist{par}, etc. In this way, the user can define state-dependent
point delays with arbitrary levels of nesting. The delay function for
the example \eqref{example_sys2} is given in Listing~\ref{sd_tau}.

\paragraph{Note: order of delays} The order of the delays
requested in \blist{sys_tau} corresponds to the order in which they
appear in \blist{xx} as passed to the functions \blist{sys_rhs} and
\blist{sys_deri}.
\paragraph{Note: difference to section~\ref{sec:consttau}} When calling
  \blist{sys_tau} for a constant delay, the value of the delay is
  returned.  This is in contrast with the definition of \blist{sys_tau}
  in section \ref{sec:consttau}, where the position in the parameter list
  is returned.

\begin{lstlisting}[frame=lines,float,label=sd_tau, caption={Listing of
  \blist{'sys_tau'} function for \eqref{example_sys2}, here called
  \blist{sd_tau}.}]
function tau=sd_tau(ind,xx,par)
if ind==1
  tau=par(10);
elseif ind==2
  tau=par(11);
elseif ind==3
  tau=2+par(5)*par(10)*xx(2,1)*xx(2,2);
elseif ind==4
 tau=1-1/(1+xx(2,3)*xx(1,1));
elseif ind==5
 tau=xx(4,1);
elseif ind==6
 tau=xx(5,1);
end
end
\end{lstlisting}

\subsubsection{Jacobian of right-hand side --- \blist{sys_deri}
  (optional, but recommended)}
\label{sec:sdjac}
The definition and functionality of this routine is equivalent to the
one described in section \ref{sec:constjac}.  We do not present here
the routine since it is quite long, see the Matlab code
\file{sd\_deri.m} in the demo example \file{sd\_demo}. If the user does
not provide a function for the Jacobians the finite-difference
approximation (\file{df\_deriv.m}) will be used by default.  However,
as for constant delays, it is recommended to provide at least
the first order derivatives with respect to the state variables using
analytical formulas.

\subsubsection{Jacobians of delays --- \blist{sys_dtau} (optional, but
  recommended)}
\label{sec:sddtau}
The routine of the format
\begin{lstlisting}
  function dtau=sd_dtau(ind,xx,par,nx,np)
\end{lstlisting}
supplies derivatives of all delays with respect to the state and
parameters. Its functionality is similar to the function
\blist{sys_deri}. Inputs:
\begin{itemize}
\item \blist{ind} (integer $\geq1$) the number of the delay,
\item \blist{par} ($n\times\blist{ind}$ matrix) is the state:
  \blist{xx(:,1)}$=x(t)$, \blist{xx(:,k)}$=x(t-\tau_{k-1})$ for
  $k=2\ldots\blist{ind}$;
\item \blist{nx} (empty, one integer or two integers) index (indices) of
  \blist{xx} with respect to which the right-hand side is to be
  differentiated;
\item \blist{np} (empty or integer) whether right-hand side is to be
  differentiated with respect to parameters.
\end{itemize}
The result $\parm{dtau}$ is a scalar, vector or matrix of partial
derivatives of the delay with number \blist{ind}, which depends on the
type of derivative requested via \blist{nx} and \blist{np}, see table
\ref{table_sd}. The resulting routine is quite long, even for the small system
\eqref{example_sys2}; see Listing~\ref{sd_dtau} in
Appendix~\ref{sec:sys:deri} for a printout of the function body.

If the user does not provide a \blist{'sys_dtau'} function then the
default routine \blist{df_derit} will be used, which implements finite
difference formulas to approximate the requested derivatives (using
several calls to \blist{sys_tau}), analogously to \blist{df_deriv}. As
in the case of \blist{sys_deri}, it is recommended to provide at least
the first order derivatives with respect to the state variables using
analytical formulas.

\begin{table}[htbp]
\begin{center}
\begin{tabular}{cc@{\hspace*{5em}}l}
  \noalign{\medskip}\hline\noalign{\smallskip}\blist{length(nx)} & \blist{length(np)}  &  \blist{dtau} 
  \\\noalign{\smallskip}\hline\noalign{\medskip}
  1         & 0      & 
  $\frac{\textstyle\partial \tau_{\blist{ind}\phantom{y}}}{\textstyle\partial x^{\,\blist{nx(1)}}}\in\RR^{n}$ \\[3ex]
  0         & 1         & 
  $\frac{\textstyle\partial \tau_{\blist{ind}\phantom{y}}}{\textstyle
    \partial^{\phantom{l}} \eta_{\blist{np(1)}}} \in\RR$ \\[3ex]
  1         & 1         & 
  $\frac{\textstyle\partial^2 \tau_{\blist{ind}\phantom{y}}}{\textstyle
    \partial x^{\blist{nx(1)}}\ \partial \eta_{\blist{np(1)}}}
  \in\RR^{n}$ \\[3ex]
  2         & 0         & 
  $\frac{\textstyle\partial}{\textstyle\partial x^{\blist{nx(2)}}}
  \left[\frac{\textstyle\partial \tau_{\blist{ind}\phantom{y}}}
    {\textstyle\partial x^{\blist{nx(1)}}}\right] \in\RR^{n\times n}$
  \\\noalign{\medskip}\hline
\end{tabular}
\caption{\label{table_sd} Results of the function
  \blist{sys_dtau} depending on its input parameters \blist{nx} and
  \blist{np} (\blist{ind}$=1,\ldots,m$).}
\end{center}
\end{table}
\subsection{Extra conditions --- \blist{sys_cond}}\label{sec:syscond}
A system routine \blist{sys_cond} with a header of the type
\begin{lstlisting}
  function [res,p]=sys_cond(point)
\end{lstlisting}
can be used to add extra conditions during corrections and
continuation, see section \ref{extra_cond} for an explanation of
arguments and outputs.

\subsection{Collecting user functions into a structure --- call
  \blist{set_funcs}}
\label{sec:funcs}
The user-provided functions are passed on as an additional argument to
all routines of \DDEBIFCODE{} (similar to standard Matlab routines
such as \blist{ode45}). This was changed in \DDEBIFCODE{}\,3.0 from
previous versions. The additional argument is a structure
\blist{funcs} containing all the handles to all user-provided
functions. In order to create this structure the user is recommended
to call the function \blist{set_funcs} at the beginning of the script
performing the bifurcation analysis:
\begin{lstlisting}
  function funcs=set_funcs(...)
\end{lstlisting}
Its argument format is in the form of name-value pairs (in arbitrary
order, similar to options at the end of a call to \blist{plot}). For
the example \eqref{example_sys} of a neuron, discussed in
section~\ref{sys_def1} and in demo \texttt{neuron} (see
\demobase{}), the call to \blist{set_funcs} could look as
follows:
\begin{lstlisting}
funcs=set_funcs('sys_rhs',neuron_sys_rhs,'sys_tau',@()[5,6,7],...
                'sys_deri',@neuron_sys_deri);
\end{lstlisting}
Note that \blist{neuron_sys_rhs} is a variable (a function handle
pointing to an anonymous function defined as in
section~\ref{sec:constrhs}), and \file{neuron\_sys\_deri.m} is the
filename in which the function providing the system derivatives are
defined (see section~\ref{sec:constjac}). The delay function
\blist{'sys_tau'} is directly specified as an anonymous function in
the call to \blist{set_funcs} (not needing to be defined in a separate
file or as a separate variable). If one does wish to not provide
analytical derivatives, one may drop the \blist{'sys_deri'} pair (then
a finite-difference approximation, implemented in \blist{df_deriv}, is
used):
\begin{lstlisting}
funcs=set_funcs('sys_rhs',neuron_sys_rhs,'sys_tau',@()[5,6,7]);
\end{lstlisting}
For the sd-DDE example \eqref{example_sys2}, the call could look as follows:
\begin{lstlisting}
funcs=set_funcs('sys_rhs',@sd_rhs,'sys_tau',@sd_tau,...
    'sys_ntau',@()6,'sys_deri',@sd_deri,'sys_dtau',@sd_dtau);  
\end{lstlisting}
Possible names to be used in the argument sequence
equal the resulting field names in \blist{funcs} (see
Table~\ref{tab:funcs} in section~\ref{sec:funcs:struct} later):
\begin{itemize}
\item \blist{'sys_rhs'}: handle of the user function providing the
  right-hand side, described in sections~\ref{sec:constrhs} and
  \ref{sec:sdrhs};
\item \blist{'sys_tau'}: handle of the user function providing the
  indices of the delays in the parameter vector for DDEs with constant
  delays, described in sections~\ref{sec:consttau}, or providing the
  values of the delays for sd-DDEs as described in
  section~\ref{sec:sdtau};
\item \blist{'sys_ntau'} (relevant for sd-DDEs only): handle of the
  user function providing the number of delays in sd-DDEs as described
  in section~\ref{sec:sdtau};
\item \blist{'sys_deri'}: handle of the user function providing the
  Jacobians of right-hand side, described in
  sections~\ref{sec:constjac} and \ref{sec:sdjac};
\item \blist{'sys_dtau'} (relevant for sd-DDEs only): handle of the user
  function providing the Jacobians of delays (for sd-DDEs), described
  in section~\ref{sec:constjac};
\item \blist{'x_vectorized'} (logical, default \blist{false}): if the
  functions in \blist{'sys_rhs'}, \blist{'sys_deri'} (if provided),
  \blist{'sys_tau'} (for sd-DDEs) and \blist{'sys_dtau'} (for sd-DDEs
  if provided) can be called with $3$d arrays in their \blist{xx}
  argument. Vectorization will speed up computations for periodic
  orbits only.
\end{itemize}
An example for a necessary modification of the right-hand side to
permit vectorization is given for the neuron example in
Listing~\ref{neuron_sys_rhs_vec}. The output \blist{funcs} is a
structure containing all user-provided functions and defaults for the
Jacobians if they are not provided. This output is passed on as first
argument to all \DDEBIFCODE{} routines during bifurcation analysis.

\section{Data structures}\label{data_structures}

In this section we describe the data structures used to define the
problem, and to present individual points, stability information, branches
of points, method parameters and plotting information.


\subsection{Problem definition (functions) structure}
\label{sec:funcs:struct}
\begin{table}[htbp]
  \centering
  \begin{tabular}[t]{l@{\hspace*{1ex}}c@{\hspace*{3ex}}p{0.5\textwidth}}
      \hline\noalign{\smallskip}
      field     & content &  default   \\\hline\noalign{\smallskip}
      \blist{'sys_rhs'} \textbf{\textsf{(mandatory)}} & function handle &  function in file \file{sys\_rhs.m} if found in 
      current working folder \\[1ex]
      \blist{'sys_ntau'} & function handle &  \blist{@()0} \\[1ex]
      \blist{'sys_tau'} \textbf{\textsf{(mandatory)}} & function handle &  function in file 
      \file{sys\_tau.m} if found in current working folder\\[1ex]
      \blist{'sys_cond'} & function handle & \blist{@(p)dummy_cond}, a built-in routine that adds no conditions\\[1ex]
      \blist{'sys_deri'} & function handle & \blist{df_deriv}, coming with 
      \DDEBIFCODE{} and using finite-difference approximation\\[1ex]
      \blist{'sys_dtau'} & function handle &  \blist{df_derit}, coming with 
      \DDEBIFCODE{} and using finite-difference approximation\\[1ex]
      \blist{'x_vectorized'} & logical &  \blist{false}\\[1ex]
      (\blist{'tp_del'}) & logical & n/a (automatically
      determined from the number of arguments expected by \blist{funcs.sys_tau})\\[1ex]
      (\blist{'sys_deri_provided'}) & logical & n/a (set automatically)\\[1ex]
      (\blist{'sys_dtau_provided'}) & logical & n/a (set automatically)
      \\\noalign{\smallskip}\hline
  \end{tabular}
  \caption{\textbf{\textsf{Problem definition structure}} containing (at least) 
    the user-provided functions. 
    Fields in brackets should not normally be set or manually changed by the user. See also section~\ref{sec:funcs}.}
  \label{tab:funcs}
\end{table}
The user-provided functions described in Section~\ref{sec:system:def}
get passed on to \DDEBIFCODE's routines collected in a single argument
\blist{funcs}, a structure containing at least the fields listed in
Table~\ref{tab:funcs}. Only fields marked with \textbf{[!]} are
mandatory (for sd-DDEs \blist{'sys_ntau'} is mandatory, too). The user
does usually not have to set the fields of this structure manually,
but calls the routine \blist{set_funcs}, which returns the structure
\blist{funcs} to be passed on to other functions. The usage of
\blist{set_funcs} and the meaning of the fields of its output are
described in detail in section~\ref{sec:funcs}. See also the tutorial
demos \texttt{neuron} and \texttt{sd\_demo} (see
\demobase{} for examples of usage.

\subsection{Point structures}\label{sec:point:struct}

Table \ref{point_structures} describes the structures used to
represent a single steady state, fold, Hopf, 
periodic and 
homoclinic/heteroclinic solution point.

\begin{table}[htbp]
{\renewcommand{\blist}[1]{\mbox{\lstinline[basicstyle={\ttfamily\small}]!#1!}}
  \begin{center}
    \begin{subtable}[b]{0.32\textwidth}\centering
      \begin{tabular}[t]{l@{\hspace*{1ex}}c}\hline\noalign{\smallskip}
        field     & content           \\\hline\noalign{\smallskip}
        \blist{'kind'}      & \blist{'stst'}            \\
        \blist{'parameter'} & $\RR^{1\times p}$ \\
        \blist{'x'}         & $\RR^{n\times 1}$ \\
        \blist{'stability'} & empty\,or\,struct\\ \\ \\\hline
      \end{tabular}
      \caption{Steady state}
    \end{subtable}
    \begin{subtable}[b]{0.32\textwidth}\centering
      \begin{tabular}[t]{l@{\hspace*{1ex}}c}\hline\noalign{\smallskip}
        field     & content           \\\hline\noalign{\smallskip} 
        \blist{'kind'}      & \blist{'fold'}            \\
        \blist{'parameter'} & $\RR^{1\times p}$ \\
        \blist{'x'}         & $\RR^{n\times 1}$ \\
        \blist{'v'}         & $\RR^{n\times 1}$ \\
        \blist{'stability'} & empty\,or\,struct \\ \\\hline
      \end{tabular}
      \caption{Steady state fold}
    \end{subtable}
    \begin{subtable}[b]{0.32\textwidth}\centering
      \begin{tabular}[t]{l@{\hspace*{1ex}}c}\hline\noalign{\smallskip}
        field     & content           \\\hline \noalign{\smallskip}
        \blist{'kind'}      & \blist{'hopf'}            \\
        \blist{'parameter'} & $\RR^{1\times p}$ \\
        \blist{'x'}         & $\RR^{n\times 1}$ \\
        \blist{'v'}         & $\CC^{n\times 1}$ \\
        \blist{'omega'}     & $\RR$             \\
        \blist{'stability'} & empty\,or\,struct\\\hline
      \end{tabular}
      \caption{Steady state Hopf}
    \end{subtable}\vspace*{5ex}
    \begin{subtable}[b]{0.45\textwidth}\centering
      \begin{tabular}[t]{lc}\hline\noalign{\smallskip}
        field     & content           \\\hline \noalign{\smallskip}
       \blist{'kind'}      & \blist{'psol'}            \\
        \blist{'parameter'} & $\RR^{1\times p}$ \\
        \blist{'mesh'}      & $[0,1]^{1\times (Ld+1)}$ or empty \\
        \blist{'degree'}    & $\NN_0$           \\
        \blist{'profile'}   & $\RR^{n\times (Ld+1)}$ \\
        \blist{'period'}    & $\RR^+_0$         \\
        \blist{'stability'} & empty or struct \\ \\ \\ \\ \\ \\ \\ \\\hline
      \end{tabular}
      \caption{Periodic orbit}
    \end{subtable}
    \begin{subtable}[b]{0.45\textwidth}\centering
      \begin{tabular}[t]{lc}\hline\noalign{\smallskip}
        field     & content           \\\hline \noalign{\smallskip}
        \blist{'kind'}      & \blist{'hcli'}            \\
        \blist{'parameter'} & $\RR^{1\times p}$ \\
        \blist{'mesh'}      & $[0,1]^{1\times (Ld+1)}$ or empty \\
        \blist{'degree'}    & $\NN_0$           \\
        \blist{'profile'}   & $\RR^{n\times (Ld+1)}$ \\
        \blist{'period'}    & $\RR^+_0$         \\
        \blist{'x1'}        & $\RR^{n}$ \\
        \blist{'x2'}        & $\RR^{n}$ \\
        \blist{'lambda_v'}  & $\CC^{s_1}$\\
        \blist{'lambda_w'}  & $\CC^{s_2}$\\
        \blist{'v'}         & $\CC^{n\times s_1}$\\
        \blist{'w'}         & $\CC^{n\times s_2}$\\
        \blist{'alpha'}     & $\CC^{s_1}$\\
        \blist{'epsilon'}   & $\RR$\\\hline
      \end{tabular}
      \caption{Connecting orbit}
    \end{subtable}
  \end{center}
}
\caption{\label{point_structures} \textbf{\textsf{Point structures}}:
Field names and corresponding content for the 
point structures used to represent steady state solutions, fold and Hopf 
points, periodic solutions and connecting orbits. Here, $n$ is 
the system dimension, $p$ is the
number of parameters, $L$ is the number of intervals used to represent
the periodic solution, $d$ is the degree of the polynomial on each
interval, $s_1$ is the number of unstable modes of $x^-$ and $s_2$ is the
number of unstable modes of $x^+$.}
\end{table}

A steady state solution is represented by the parameter values
$\eta$ (which contain also the constant delay values, 
see section \ref{sec:system:def})
and $x^*$. A fold bifurcation is represented by the parameter
values $\eta$, its position $x^*$ and a null-vector of the
characteristic matrix $\Delta(0)$. A Hopf bifurcation is represented 
by the parameter
values $\eta$, its position $x^*$, a frequency $\omega$
and a (complex) null-vector of the
characteristic matrix $\Delta(\i \omega)$.

A periodic solution is represented by the parameter
values $\eta$, the period $T$ and 
a time-scaled profile $x^*(t/T)$ on a mesh in [0,1].
The mesh is an ordered 
collection of \define{interval points} $\{0=t_0<t_1<\ldots<t_L=1\}$
and \define{representation points} $t_{i+\frac{j}{d}}$, $i=0,\ldots,L-1$,
$j=1,\ldots,d-1$ which need to be chosen in function of the interval points as
\[
t_{i+\frac{j}{d}}=t_i+\frac{j}{d}(t_{i+1}-t_i).
\]
\paragraph{Warning:
this assumption is not checked but needs to be fulfilled for
correct results!}
The profile is a continuous piecewise polynomial on the mesh. 
More specifically, it is
a polynomial of degree $d$ on each
subinterval $[t_i,t_{i+1}]$, $i=0,\ldots,L-1$.
Each of these polynomials is uniquely represented
by its values at the points $\{t_{i+\frac{j}{d}}\}_{j=0,\ldots,d}$.
Hence the complete profile is represented by
its value at all the mesh points,
\[
x^*(t_{i+\frac{j}{d}}),\ i=0,\ldots,L-1,\ j=0,\ldots,d-1;
\mathrm{\ and\ } x^*(t_L).
\]
Because polynomials on adjacent intervals share the value
at the common interval point, this representation is automatically
continuous (it is, however, not continuously differentiable).
(As indicated in table \ref{point_structures}, the mesh may be empty, 
which indicates the use of an equidistant, fixed mesh.)

A connecting orbit is represented by the parameter values $\eta$, the
period $T$, a time-scaled profile $x^*(t/T)$ on a mesh in [0,1], the
steady states $x^-$ and $x^+$ (fields \blist{'x1'} and \blist{'x2'} in
the data structure), the unstable eigenvalues of these steady states,
$\lambda^-$ and $\lambda^+$ (fields \blist{'lambda_v'} and
\blist{'lambda_w'} in the data structure), the unstable right
eigenvectors of $x^-$ (\blist{'v'}), the unstable left eigenvectors of
$x^+$ (\blist{'w'}), the direction in which the profile leaves the unstable
manifold, determined by $\alpha$, and the distance of the first point
of the profile to $x^-$, determined by $\epsilon$.  For the mesh and
profile, the same remarks as in the case of periodic solutions hold.

The point structures are used as input to the point
manipulation routines (layer 2) and are used inside
the branch structure (see further). The order of the fields in the point 
structures is important (because they are used as elements
of an array inside the branch structure).
No such restriction holds for the other structures (method, plot and
branch) described in the rest of this section.

\subsection{Stability structures}\label{sec:stab:struct}

Most of the point structures contain a 
field \blist{'stability'} storing eigenvalues or Floquet multipliers. (The exception is  
the \blist{'hcli'} structure for which stability does not really make
sense.) During bifurcation analysis the computation of stability is
typically performed as a separate step, after computation of the
solution branches, because stability computation can easily be more
expensive than the solution finding. If no stability has been computed
yet, the field \blist{'stability'} is empty, otherwise, it contains
the computed stability information in the form described in
Table~\ref{stab_structures}.

\begin{table}[htbp]
  \begin{center}
    \begin{subtable}[t]{0.45\textwidth}
      \begin{tabular}[t]{lc}\hline\noalign{\smallskip}
        field & content     \\\hline\noalign{\smallskip}
        \blist{'h'}     & $\RR$       \\
        \blist{'l0'}    & $\CC^{n_l}$ \\
        \blist{'l1'}    & $\CC^{n_c}$ \\
        \blist{'n1'}    & $(\{-1\}\cup\NN_0)^{n_c}$ or empty\\\hline
      \end{tabular}
      \caption{Structure in field \blist{'stability'} for steady
        state, fold and Hopf points of Table~\ref{point_structures}}
    \end{subtable}\qquad
    \begin{subtable}[t]{0.45\textwidth}
      \begin{tabular}[t]{lc}\hline\noalign{\smallskip}
        field & content     \\\hline\noalign{\smallskip}
        \blist{'mu'}    & $\CC^{n_m}$ \\ \\ \\ \\\hline
      \end{tabular}
      \caption{Structure in field \blist{'stability'} for periodic
        orbit points of Table~\ref{point_structures}}
    \end{subtable}
  \end{center}
  \caption{\label{stab_structures}
    \textbf{\textsf{Stability structures}} for roots of the characteristic equation
    (in steady state, fold and Hopf structures) (left)
    and for Floquet multipliers (in the periodic solutions structure) (right). 
    Here, $n_l$ is the number of approximated roots,
    $n_c$ is the number of corrected roots and $n_m$ is the number
    of Floquet multipliers.}
\end{table}
For steady state, fold and Hopf points, approximations to the
rightmost roots of the characteristic equation are provided in field
\blist{'l0'} in order of decreasing real part.  The steplength that
was used to obtain the approximations is provided in field
\blist{'h'}. Corrected roots are provided in field \blist{'l1'} and
the number of Newton iterations applied for each corrected root in a
corresponding field \blist{'n1'}.  If unconverged roots are discarded,
\blist{'n1'} is empty and the roots in \blist{'l1'} are ordered with
respect to real part; otherwise the order in \blist{'l1'} corresponds
to the order in \blist{'l0'} and an element $-1$ in \blist{'n1'}
signals that no convergence was reached for the corresponding root in
\blist{'l0'} and the last computed iterate is stored in \blist{'l1'}.
The collection of uncorrected roots presents more accurate yet less
robust information than the collection of approximate roots, see
section \ref{code_num_methods}. For periodic solutions only
approximations to the Floquet multipliers are provided in a field
\blist{'mu'} (in order of decreasing modulus). As the characteristic
matrix is not analytically available, \DDEBIFCODE{} does not offer an
additional correction.

\subsection{Method parameters}\label{sec:method:struct}
\begin{table}[htbp]
\begin{center}
\begin{tabular}{l@{\hspace*{2em}}c@{\hspace*{2em}}c}\hline\noalign{\smallskip}
field                      & content     & default value  \\\hline\noalign{\smallskip}
\blist{'newton_max_iterations'}    & $\NN_0$     & \blist{5}, \blist{5}, \blist{5}, \blist{5}, \blist{10} \\
\blist{'newton_nmon_iterations'}   & $\NN_{\phantom{0}}$       & \blist{1} \\
\blist{'halting_accuracy'}         & $\RR^+$     & \blist{1e-10}, \blist{1e-9}, \blist{1e-9}, \blist{1e-8}, \blist{1e-8} \\[1ex]
\blist{'minimal_accuracy'}         & $\RR^+_0$   & \texttt{\ }\blist{1e-8}, \blist{1e-7}, \blist{1e-7}, \blist{1e-6}, \blist{1e-6} \\
\blist{'extra_condition'}          & $\{0,1\}$   & \blist{0} \\
\blist{'print_residual_info'}      & $\{0,1\}$   & \blist{0}\\[2ex]

$^*$\blist{'phase_condition'}             & $\{0,1\}$          & \blist{1} \\
$^*$\blist{'collocation_parameters'}      & $[0,1]^d$ or empty & empty \\
$^*$\blist{'adapt_mesh_before_correct'} & $\NN$              & \blist{0} \\
$^*$\blist{'adapt_mesh_after_correct'}  & $\NN$              & \blist{3} 
\\\noalign{\smallskip}\hline
\end{tabular}
\end{center}
\caption{\label{point_method_structures}
  \textbf{\textsf{Point method structure}}: fields and possible values. 
  When different,
  default values are given in the order \blist{'stst'},\blist{'fold'},\blist{'hopf'},\blist{'psol'}, \blist{'hcli'}. Fields marked with and asterisk ($^*$) are needed and present for points of type \blist{'psol'} and \blist{'hcli'} only.}
\end{table} 
To compute a single steady state, fold, Hopf, periodic 
or connecting
orbit solution point,
several method parameters have to be passed to the appropriate routines.
These parameters are collected into a structure with the fields
given in Table \ref{point_method_structures}.
\begin{table}[htbp]
\begin{center}
\begin{tabular}{l@{\hspace*{1em}}c@{\hspace*{1em}}c}\hline\noalign{\smallskip}
 field                        & content              & default value  \\\hline\noalign{\smallskip} 
\textbf{For steady state, fold and Hopf}\\\noalign{\smallskip}
\blist{'lms_parameter_alpha'}        & $\RR^k$              & \blist{time_lms('bdf',4)} \\
\blist{'lms_parameter_beta'}         & $\RR^k$              & \blist{time_lms('bdf',4)} \\
\blist{'lms_parameter_rho'}          & $\RR^+_0$            & \blist{time_saf(alpha,beta,0.01,0.01)} \\
\blist{'interpolation_order'}         & $\NN_0$              & \blist{4} \\
\blist{'minimal_time_step'}          & $\RR^+_0$            & \blist{0.01} \\
\blist{'maximal_time_step'}          & $\RR^+_0$            & \blist{0.1} \\
\blist{'max_number_of_eigenvalues'} & $\NN_0$              & \blist{100} \\
\blist{'minimal_real_part'}          & $\RR$ or empty       & empty \\
\blist{'max_newton_iterations'}      & $\NN$                & \blist{6} \\ 
\blist{'root_accuracy'}               & $\RR^+_0$            & \blist{1e-6} \\
\blist{'remove_unconverged_roots'}   & $\{0,1\}$            & \blist{1} \\
\blist{'delay_accuracy'}              & $\RR_0^-$            & \blist{-1e-8}
\\\noalign{\medskip}
\textbf{For periodic orbit}\\\noalign{\smallskip}
\blist{'collocation_parameters'}      & $[0,1]^d$ or empty   & empty \\
\blist{'max_number_of_eigenvalues'} & $\NN$                & \blist{100} \\
\blist{'minimal_modulus'}             & $\RR^+$              & \blist{0.01} \\
\blist{'delay_accuracy'}              & $\RR_0^-$            & \blist{-1e-8}
\\\noalign{\smallskip}\hline
\end{tabular}
\end{center}
\caption{\label{meth_stab_struct}
\textbf{\textsf{Stability method structures}}: fields and possible values
for the approximation and correction of roots of the characteristic
equation (top), or
for the approximation Floquet multipliers (bottom).
The LMS-parameters are default set to the fourth order backwards
differentiation LMS-method.
The last row in both parts is only used for sd-DDEs.}
\end{table}

For the computation of periodic solutions, additional fields are
necessary, marked with and asterisk ($^*$) in Table \ref{point_method_structures}.
The meaning of the different fields in Table
\ref{point_method_structures} is explained in section
\ref{code_num_methods}.

\begin{table}
\begin{center}
\begin{tabular}{l@{\hspace*{2em}}c@{\hspace*{2em}}c}\hline\noalign{\smallskip}
field                      & content         & default value  \\\hline\noalign{\smallskip} 
\blist{'steplength_condition'}      & $\{0,1\}$       & \blist{1}     \\
\blist{'plot'}                       & $\{0,1\}$       & \blist{1}     \\
\blist{'prediction'}                 & $\{1\}$         & \blist{}1     \\
\blist{'steplength_growth_factor'} & $\RR^+_0$       & \blist{1.2}   \\
\blist{'plot_progress'}             & $\{0,1\}$       & \blist{1}     \\
\blist{'plot_measure'}              & struct or empty & empty \\
\blist{'halt_before_reject'}       & $\{0,1\}$       & \blist{0}
\\\noalign{\smallskip}\hline
\end{tabular}
\end{center}
\caption{\label{continuation_structure}
\textbf{\textsf{Continuation method structure}}: fields and possible values.}
\end{table}
Parameters controlling the pseudo-arclength continuation (using secant
approximations for tangents) are stored in a structure of the form
given in Table~\ref{continuation_structure}.
Similarly, for the approximation and correction of roots of the 
characteristic equation respectively for the computation of the
Floquet multipliers
method parameters are passed using a structure of the form given
in table \ref{meth_stab_struct}.
\subsection{Branch structures}\label{sec:branch:struct}
\begin{table}
\begin{center}
\begin{tabular}{l@{\hspace*{2em}}l@{\hspace*{2em}}cl}\hline\noalign{\smallskip}
field     & subfield     & content                    \\\hline \noalign{\smallskip}
\blist{'point'}     &              & array of points &(s. Table~\ref{point_structures})      \\[1ex]
\blist{'method'}    & \blist{'point'}        & point method struct &(s. Table~\ref{point_method_structures})        \\
        & \blist{'stability'}    & stability method struct  &(s. Table~\ref{meth_stab_struct})   \\ 
       & \blist{'continuation'} & continuation method struct &(s. Table~\ref{continuation_structure}) \\[1ex] 
\blist{'parameter'} & \blist{'free'}         & $\NN^{p_f}$                \\ 
 & \blist{'min_bound'}   & $[\NN\ \RR]^{p_i}$        \\
 & \blist{'max_bound'}   & $[\NN\ \RR]^{p_a}$        \\
 & \blist{'max_step'}    & $[\NN\ \RR]^{p_s}$        \\
\noalign{\smallskip}\hline
\end{tabular}\end{center}
\caption{\label{branch_struct}\textbf{\textsf{Branch structure}}: fields and 
possible values. Here, $p_f$ is the number of free parameters; 
$p_i$, $p_a$ and $p_s$ are the number of minimal parameter values, 
maximal parameter values respectively 
maximal parameter steplength values. If any of these values are zero,
the corresponding subfield is empty.}
\end{table}
A branch consists of an ordered array of points (all of the same type), 
and three method structures
containing point method parameters, continuation parameters
respectively stability computation parameters,
see table \ref{branch_struct}.

The branch structure has three fields. One, called \blist{'point'},
which contains an array of point structures, one, called
\blist{'method'}, which is itself a structure containing three
subfields and a third, called \blist{'parameter'} which contains four
subfields.  The three subfields of the method field are again
structures. The first, called \blist{'point'}, contains point method
parameters as described in Table \ref{point_method_structures}.  The
second, called \blist{'stability'}, contains stability method
parameters as described in Table \ref{meth_stab_struct} and the third,
called \blist{'continuation'}, contains continuation method parameters
as described in Table \ref{continuation_structure}.  Hence the branch
structure incorporates all necessary method parameters which are thus
automatically kept when saving a branch variable to file.  The
parameter field contains a list of free parameter numbers which are
allowed to vary during computations, and a list of parameter bounds
and maximal steplengths. Each row of the bound and steplength
subfields consists of a parameter number (first element) and the value
for the bound or steplength limitation. Examples are given in demo
\texttt{neuron} (see \demobase{}).

A default, empty branch structure can be obtained by passing
a list of free parameters and the point kind 
(as \blist{'stst'}, \blist{'fold'}, \blist{'hopf'}, \blist{'psol'}
or \blist{'hcli'})
to the function \blist{df_brnch}. A minimal bound zero is then set
for each constant delay if the function \blist{sys_tau} is defined
as in section~\ref{sys_def1} (i.e.~for DDEs). The method contains 
default parameters
(containing appropriate point, stability and continuation fields)
obtained from the function \blist{df_mthod} with as only argument the type
of solution point.

\subsection{Scalar measure structure}\label{sec:meas:struct}

After a branch has been computed some possibilities are offered to
plot its content. For this a (scalar) measure structure is used which
defines what information should be taken and how it should be
processed to obtain a measure of a given point (such as the amplitude
of the profile of a periodic solution, etc\ldots); see Table
\ref{measure_structure}.  The result applied to a variable
\blist{point} is to be interpreted as
\begin{lstlisting}[frame=none]
scalar_measure=func(point.field.subfield(row,col));
\end{lstlisting}
where \blist{'field'} presents the field to select,
\blist{'subfield'} is empty or presents the subfield to select, 'row'
presents the row number or contains one of the functions mentioned in
table \ref{measure_structure}.  These functions are applied columnwise
over all rows.  The function \blist{'all'} specifies that the all
rows should be returned.  The meaning of \blist{'col'} is similar to
\blist{'row'} but for columns.  To avoid ambiguity it is required
that either \blist{'row'} or \blist{'col'} contains a number or that
both contain the function \blist{'all'}.  If nonempty, the function
'func' is applied to the result.  Note that \blist{'func'} can be a
standard Matlab function as well as a user written function. Note also
that, when using the value \blist{'all'} in the fields \blist{'col'}
and/or \blist{'row'} it is possible to return a non-scalar measure
(possibly but not necessarily further processed by \blist{'func'}).

\begin{table}
\begin{center}
\begin{tabular}{l@{\hspace*{0em}}c@{\hspace*{1em}}c}\hline\noalign{\smallskip}
field  & content & meaning      \\\hline \noalign{\smallskip}
\blist{'field'} & \{\blist{'parameter'},\blist{'x'},\blist{'v'},\blist{'omega'},\,\ldots & first field to select\\     
       &\qquad\blist{'profile'},\blist{'period'},\blist{'stability'} \ldots\} &  from a point struct\\
\blist{'subfield'} & \{'\ ',\blist{'l0'},\blist{'l1'},\blist{'mu'}\} & empty string or 2nd field to select \\
\blist{'row'}     & $\NN$ or \{\blist{'min'},\blist{'max'},\blist{'mean'},\blist{'ampl'},\blist{'all'}\} & row index \\
\blist{'col'}     & $\NN$ or \{\blist{'min'},\blist{'max'},\blist{'mean'},\blist{'ampl'},\blist{'all'}\} & column index \\
\blist{'func'}   & \{\blist{''},\blist{'real'},\blist{'imag'},\blist{'abs'}\} & function to apply\\\noalign{\smallskip}\hline
\end{tabular}
\end{center}
\caption{\label{measure_structure}
  \textbf{\textsf{Measure structure}}: fields, content and meaning of
  a structure describing
  a measure of a point.}
\end{table}

\section{Point manipulation}\label{point_manipulation}

Several of the point manipulation routines have already been used in
the previous section.  Here we outline their functionality and input
and output parameters.  A brief description of parameters is also
contained within the source code and can be obtained in Matlab using
the $\parm{help}$ command. Note that a vector of zero elements
corresponds to an empty matrix (written in Matlab as $[]$). 
\begin{lstlisting}
function [point,success]=p_correc(...
          funcs,point0,free_par,step_cnd,method,adapt,previous,d_nr,tz)  
\end{lstlisting}
\noindent Function \blist{p_correc} corrects a given point.
\begin{itemize}
\item \blist{funcs}: structure of user-defined functions, defining the
  problem (created, for example, using \blist{set_funcs}).
\item \blist{point0}: initial, approximate solution point as a point
  structure (see table \ref{point_structures}).
\item \blist{free_par}: a vector of zero, one or more free parameters.
\item \blist{step_cnd}: a vector of zero, one or more linear
  steplength conditions. Each steplength condition is assumed
  fulfilled for the initial point and hence only the coefficients of
  the condition with respect to all unknowns are needed. These
  coefficients are passed as a point structure (see table
  \ref{point_structures}).  This means that for, e.g., a steady state
  solution point \blist{p} the $i$-th steplength condition enforces that
  \begin{lstlisting}
    step_cnd(i).parameter*(p.parameter-point0.parameter)'+...
    step_cnd(i).x'*(p.x-point0.x)
  \end{lstlisting}
  is zero. Similar formulas hold for the other solution types.
\item \blist{method}: a point method structure containing the method
  parameters (see table \ref{point_method_structures}).
\item \blist{adapt} (optional): if zero or absent, do not use adaptive
  mesh selection (for periodic solutions); if one, correct, use
  adaptive mesh selection and recorrect.
\item \blist{previous} (optional): for periodic solutions and
  connecting orbits: if present and not empty, minimize phase shift
  with respect to this point. Note that this argument should always be
  present when correcting solutions for sd-DDEs, since in that case
  the argument \blist{d_nr} always needs to be specified.  In the case
  of steady state, fold or Hopf-like points, one can just enter an
  empty vector.
\item \blist{d\_nr}: (only for equations with state-dependent delays)
  if present, number of a negative state-dependent delay.
\item \blist{tz}: (only for equations with state-dependent delays and
  periodic solutions) if present, a periodic solution is computed such
  that $\tau_{\blist{tz}}=0$ and ${\d}\tau_{\blist{tz}}/{\d}t=0$, where
  $\tau$ is a negative state-dependent delay with number
  \blist{d_nr}. For steady state solutions, a solution corresponding
  to $\tau=0$ is computed.
\item \blist{point}: the result of correcting \blist{point0} using the
  method parameters, steplength condition(s) and free parameter(s)
  given. Stability information present in \blist{point0} is not passed
  onto \blist{point}.  If divergence occurred, \blist{point} contains
  the final iterate.
\item \blist{success}: nonzero if convergence was detected (that is,
  if the requested accuracy has been reached).
\end{itemize}
\begin{lstlisting}
function stability=p_stabil(funcs,point,method)  
\end{lstlisting}
\noindent Function \blist{p_stabil} computes stability of a given
point by approximating its stability-determining eigenvalues.
\begin{itemize}
\item \blist{funcs}: structure of user-defined functions, defining the
  problem (created, for example, using \blist{set_funcs}).
\item \blist{point}: a solution point as a point structure (see table
  \ref{point_structures}).
\item \blist{method}: a stability method structure (see table
  \ref{meth_stab_struct}).
\item \blist{stability}: the computed stability of the point through a
  collection of approximated eigenvalues (as a structure described in
  table \ref{stab_structures}).  For steady state, fold and Hopf
  points both approximations and corrections to the rightmost roots of
  the characteristic equation are provided.  For periodic solutions
  approximations to the dominant Floquet multipliers are computed.
\end{itemize}
\begin{lstlisting}
function p_splot(point)  
\end{lstlisting}
\noindent Function \blist{p_splot} plots the characteristic roots
respectively Floquet multipliers of a given point (which should
contain nonempty stability information).  Characteristic root
approximations and Floquet multipliers are plotted using '$\times$',
corrected characteristic roots using '$*$'.

\begin{lstlisting}
function stst_point=p_tostst(funcs,point)
function fold_point=p_tofold(funcs,point)
function hopf_point=p_tohopf(funcs,point,excludefreqs)
function [psol_point,stepcond]=p_topsol(funcs,point,ampl,degree,nr_int)
function [psol_point,stepcond]=p_topsol(funcs,point,ampl,coll_points)
function [psol_point,stepcond]=p_topsol(funcs,hcli_point)
function hcli_point=p_tohcli(funcs,point)  
\end{lstlisting}
\noindent The functions \blist{p_tostst}, \blist{p_tofold},
\blist{p_tohopf}, \blist{p_topsol} and \blist{p_tohcli} convert a
given point into an approximation of a new point of the kind indicated
by their name. They are used to switch from a steady state point to a
Hopf point or fold point, from a Hopf point to a fold point or vice
versa, from a (nearby double) Hopf point to the second Hopf point,
from a Hopf point to the emanating branch of periodic solutions, from
a periodic solution near a period doubling bifurcation to the
period-doubled branch and from a periodic solution near a homoclinic
orbit to this homoclinic orbit.  The function \blist{p_tostst} is also
capable of extracting the initial and final steady states from a
connecting orbit. 

The additional argument \blist{excludefreqs} of function
\blist{p_tohopf} controls which Hopf frequency is chosen if several
are possible. Function \blist{p_tohopf} calls \blist{p_correc} with an
initial guess for the Hopf freqency equal to the eigenvalue closest to
the imaginary axis. For each entry in the vector \blist{excludefreqs}
(of non-negative real numbers) the eigenvalue with imaginary part
closest to this entry will be removed from consideration. This becomes
useful in systems with large delays. Then equilibria have a large
number of eigenvalues close to the imaginary axis, and
correspondingly, the system experiences many Hopf bifurcations over
short parameter intervals. These Hopf bifurcations can be picked up by
calling \blist{p_tohopf} repeatedly, and excluding frequencies one
after another.

When starting a periodic solution branch from a Hopf point, an
equidistant mesh is produced with \blist{nr_int} intervals and
piecewise polynomials of degree \blist{col_degree} and a steplength
condition \blist{stepcond} is returned which should be used (together
with a corresponding free parameter) in correcting the returned
point. This steplength condition (normally) prevents convergence back
to the steady state solution (as a degenerate periodic solution of
amplitude zero). When jumping to a period-doubled branch, a
period-doubled solution profile is produced using \blist{coll_points}
for collocation points and a mesh which is the (scaled) concatenation
of two times the original mesh.  A steplength condition is returned
which (normally) prevents convergence back to the single period
branch.

When jumping from a homoclinic orbit to a periodic solution, the
steplength condition prevents divergence, by keeping the period fixed.
When extracting the steady states from a connecting orbit, an array is
returned in which the first element is the initial steady state, and
the second element is the final steady state.

\begin{lstlisting}
function rm_point=p_remesh(point,new_degree,new_mesh)  
\end{lstlisting}
\noindent Function \blist{p_remesh} changes the piecewise polynomial
representation of a given periodic solution point.
\begin{itemize}
\item \blist{point}: initial point, containing old mesh, old degree
  and old profile.
\item \blist{new_degree}: new degree of piecewise polynomials.
\item \blist{new_mesh}: mesh for new representation of periodic
  solution profile either as a (non-scalar) row vector of mesh points
  (both interval and representation points, with the latter chosen
  equidistant between the former, see section \ref{data_structures})
  or as the new number of intervals.  In the latter case the new mesh
  is adaptively chosen based on the old profile.
\item \blist{rm_point}: returned point containing new degree, new mesh
  and an appropriately interpolated (but uncorrected!) profile.
\end{itemize} 
\begin{lstlisting}
function tau_eva=p_tau(funcs,point,d_nr,t)  
\end{lstlisting}
\noindent Function \blist{p_tau} evaluates state-dependent delay(s) 
with number(s) \blist{d_nr}. 
\begin{itemize}
\item \blist{funcs}: structure of user-defined functions, defining the
  problem (created, for example, using \blist{set_funcs}).
\item \blist{point}: a solution point as a point structure.
\item \blist{d_nr}: number(s) of delay(s) (in increasing order) 
to evaluate.
\item \blist{t} (absent for steady state solutions and optional for
  periodic solutions): mesh (a time point or a number of time
  points). If present, delay function(s) are evaluated at the points
  of \blist{t}, otherwise at the \blist{point.mesh} (if
  \blist{point.mesh} is empty, an equidistant mesh is used).
\item \blist{tau_eva}: evaluated values of delays (at \blist{t}). 
\end{itemize}
The following routines are used within branch routines but
are less interesting for the general user.

\begin{lstlisting}
function sc_measure=p_measur(p,measure)  
\end{lstlisting}
\noindent Function \blist{p_measur} computes the (scalar) measure
\blist{measure} of the given point \blist{p} (see table
\ref{measure_structure}).

\begin{lstlisting}
function p=p_axpy(a,x,y)  
\end{lstlisting}
\noindent Function \blist{p_axpy} performs the axpy-operation on
points. That is, it computes \blist{p}=\blist{a}\blist{x}+\blist{y} where 
\blist{a} is a scalar, and
\blist{x} and \blist{y} are two point structures of the same type.
\blist{p} is the result of the operation on all appropriate
fields of the given points.
If \blist{x} and \blist{y} are 
solutions on different meshes, interpolation
is used and the result is obtained on the mesh of \blist{x}.
Stability information, if present, is not passed onto \blist{p}.

\begin{lstlisting}
function n=p_norm(point)  
\end{lstlisting}
\noindent Function \blist{p_norm} computes some 
norm of a given point structure.

\begin{lstlisting}
function normalized_p=p_normlz(p)  
\end{lstlisting}
\noindent Function \blist{p_normlz} performs some normalization on the
given point structure \blist{p}. In particular, fold, Hopf and
connecting orbit determining eigenvectors are scaled to norm 1.

\begin{lstlisting}
function [delay_nr,tz]=p_tsgn(point)  
\end{lstlisting}
\noindent Function \blist{p_tsgn} detects a first negative state-dependent 
delay.
\begin{itemize}
\item \blist{point}: a solution point as a point structure.
\item \blist{delay_nr}: number of the first (and only the first !) 
detected negative delay $\tau$.
\item \blist{tz} (only for periodic solutions): $\blist{tz}\in [0,1]$ 
is a (time) point such that the delay function $\tau(t)$ has its 
minimal value near this point. To compute \blist{tz}, a refined mesh 
is used in the neighbourhood of the minimum of the delay function.
This point is later used to compute a periodic solution
such that $\tau_{\blist{tz}}=0$ and ${\d}\tau_{\blist{tz}}/{\d}t=0$. 
\end{itemize}

\section{Branch manipulation}\label{branch_manipulation}

Usage of most of the branch manipulation routines is illustrated in
the demos \texttt{neuron} and \texttt{sd\_demo} (see
\demobase{}).  Here we outline their functionality and
input and output variables.  As for all routines in the package, a
brief description of the parameters is also contained within the
source code and can be obtained in Matlab using the \blist{help}
command.
\begin{lstlisting}
function [c_branch,succ,fail,rjct]=br_contn(funcs,branch,max_tries)  
\end{lstlisting}
\noindent The function \blist{br_contn} computes (or rather 
extends) a branch of solution
points. 
\begin{itemize}
\item \blist{funcs}: structure of user-defined functions, defining the
  problem (created, for example, using \blist{set_funcs}).
\item \blist{branch}: initial branch 
containing at least two points and computation, stability and 
continuation method parameter
structures and a free parameter structure as described in 
table \ref{branch_struct}. 
\item \blist{max_tries}:
maximum number of corrections allowed.
\item \blist{c_branch}:
the branch returned contains a copy of the initial branch plus
the extra points computed (starting from the end of the point array in the
initial branch). 
\item \blist{succ}: number of successful corrections.
\item \blist{fail}: number of failed corrections.
\item \blist{rjct}: number of rejected points.
\end{itemize}
Note also that successfully computed points are normalized using the procedure
\blist{p_normlz} (see section \ref{point_manipulation}). 
\begin{lstlisting}
function br_plot(branch,x_measure,y_measure,line_type)  
\end{lstlisting}
\noindent Function \blist{br_plot} plots a branch (in the current figure). 
\begin{itemize}
\item \blist{branch}: branch to plot (see table \ref{branch_struct}).
\item \blist{x_measure}: (scalar) measure to produce plotting quantities
for the x-axis (see table \ref{measure_structure}). 
If empty, the point number is used to plot against.
\item \blist{y_measure}: (scalar) measure to produce plotting quantities
for the y-axis (see table \ref{measure_structure}). 
If empty, the point number is used to plot against.
\item \blist{line_type} (optional): line type to plot with.
\end{itemize}

\begin{lstlisting}
function [x_measure,y_measure]=df_measr(stability,branch)
function [x_measure,y_measure]=df_measr(stability,par_list,kind)  
\end{lstlisting}
\noindent Function \blist{df_measr} returns default measures for
plotting.
\begin{itemize}
\item \blist{stability}: nonzero if measures are required to plot
  stability information.
\item \blist{branch}: a given branch (see table \ref{branch_struct})
  for which default measures should be constructed.
\item \blist{par_list}: a list of parameters for which default
  measures should be constructed.
\item \blist{kind}: a point type for which default measures should be
  constructed.
\item \blist{x_measure}: default scalar measure to use for the x-axis.
  \blist{x_measure} is chosen as the first parameter which varies
  along the branch or as the first parameter of \blist{par_list}.
\item \blist{y_measure}: default scalar measure to use for the y-axis.
  If \blist{stability} is zero, the following choices are made for
  \blist{y_measure}. For steady state solutions, the first component
  which varies along the branch; for fold and Hopf bifurcations the
  first parameter value (different from the one used for
  \blist{x_measure}) which varies along the branch. For periodic
  solutions, the amplitude of the fist varying component.  If
  \blist{stability} is nonzero, \blist{y_measure} selects the real
  part of the characteristic roots (for steady state solutions, fold
  and Hopf bifurcations) or the modulus of the Floquet multipliers
  (for periodic solutions).
\end{itemize}

\begin{lstlisting}
function st_branch=br_stabl(funcs,branch,skip,recompute)  
\end{lstlisting}
\noindent Function \blist{br_stabl} computes stability information
along a previously computed branch. 
\begin{itemize}
\item \blist{funcs}: structure of user-defined functions, defining the
  problem (created, for example, using \blist{set_funcs}).
\item \blist{branch}: given branch (see table \ref{branch_struct}).
\item \blist{skip}: number of points to skip between stability
  computations.  That is, computations are performed and stability
  field is filled in every \blist{skip}+1-th point.
\item \blist{recompute}: if zero, do not recompute stability
  information present. If nonzero, discard and recompute old stability
  information present (for points which were not skipped).
\item \blist{st_branch}: a copy of the given branch whose
  (non-skipped) points contain a non-empty stability field with
  computed stability information (using the method parameters
  contained in \blist{branch}).
\end{itemize}

\begin{lstlisting}
function t_branch=br_rvers(branch)  
\end{lstlisting}
\noindent To continue a branch in the other direction (from the
beginning instead of from the end of its point array),
\blist{br_rvers} reverses the order of the points in the branches
point array.

\begin{lstlisting} 
function recmp_branch=br_recmp(funcs,branch,point_numbers) 
\end{lstlisting}
\noindent Function \blist{br_recmp} recomputes part of a branch.
\begin{itemize}
\item \blist{funcs}: structure of user-defined functions, defining the
  problem (created, for example, using \blist{set_funcs}).
\item \blist{branch}: initial branch (see table \ref{branch_struct}).
\item \blist{point_numbers} (optional): vector of one or more point
  numbers which should be recomputed. Empty or absent if the complete
  point array should be recomputed.
\item \blist{recmp_branch}: a copy of the initial branch with points
  who were (successfully) recomputed replaced. If a recomputation
  fails, a warning message is given and the old value remains present.
\end{itemize}
This routine can, e.g., be used after changing some method parameters
within the branch method structures.

\begin{lstlisting}
function [col,lengths]=br_measr(branch,measure)  
\end{lstlisting}
\noindent Function \blist{br_selec} computes a measure along a branch.
\begin{itemize}
\item \blist{branch}: given branch (see table \ref{branch_struct}).
\item \blist{measure}: given measure (see table \ref{measure_structure}).
\item \blist{col}: the collection of measures taken along the
branch (over its point array) ordered row-wise. Thus, a column vector
is returned if \blist{measure} is scalar. Otherwise,
\blist{col} contains a matrix.
\item \blist{lengths}: vector of lengths of the measures along the
  branch.  If the measure is not scalar, it is possible that its
  length varies along the branch (e.g.~when plotting rightmost
  characteristic roots). In this situation \blist{col} is a matrix
  with number of columns equal to the maximal length of the measures
  encountered.  Extra elements of \blist{col} are automatically put to
  zero by Matlab.  \blist{lengths} can then be used to prevent
  plotting of extra zeros.
\end{itemize}

\section{Numerical methods}\label{numerical_methods}\label{code_num_methods}
This section contains short descriptions of the numerical methods 
for DDEs and the method parameters used in {\DDEBIFCODE}. 
More details on the methods can be found in the
articles \cite{Luzy96,Enge99a,Enge99b,en_d01,engel01,homoclinic} 
or in \cite{Enge00}. For details on applying these methods to bifurcation
analysis of sd-DDEs see \cite{luz01}.

\subsection{Determining systems}\label{determining_systems}

Below we state the determining systems used to compute and
continue steady state solutions, steady state fold and Hopf 
bifurcations, periodic solutions and connecting orbits of systems of delay
differential equations.

For each determining system we mention the number of free 
parameters necessary to obtain (generically) isolated 
solutions. 
In the package,
the necessary number of free parameters
is further raised by the number of
steplength conditions plus the number of extra conditions used.
This choice ensures 
the use of square Jacobians during Newton iteration. 
If, on the other hand, the number of free parameters, 
steplength conditions and extra conditions
are not appropriately matched Newton iteration solves systems with a   
non-square Jacobian (for which Matlab uses an
over- or under-determined
least squares procedure). 
If possible, it is better to avoid such a situation.

\paragraph{Steady state solutions}
A steady state solution $x^*\in\RR^n$ is determined from the following
$n$-dimensional determining system with no free parameters.
\begin{equation}\label{determ_stst}
f(x^*,x^*,\ldots,x^*,\eta)=0.
\end{equation}

\paragraph{Steady state fold bifurcations}
Fold bifurcations, $(x^*\in\RR^n,v\in\RR^n)$ are determined 
from the following
$2n+1$-dimensional determining system using one free
parameter.
\begin{equation}\label{determ_fold}
\begin{aligned}
0&=f(x^*,x^*,\ldots,x^*,\eta) \\
0&=\Delta(x^*,\eta,0)v\\
0&=c^\T v-1
\end{aligned}
\end{equation}
(see \eqref{eq:deltadef} for the definition of the characteristic
matrix $\Delta$). Here, $c^\T v-1=0$ presents a suitable normalization
of $v$.  The vector $c\in\RR^n$ is chosen as $c=v^{(0)}/({v^{(0)}}^\T
v^{(0)})$, where $v^{(0)}$ is the initial value of $v$.

\paragraph{Steady state Hopf bifurcations}
Hopf bifurcations, $(x^*\in\RR^n,v\in\CC^n,\omega\in\RR)$ are
determined from the following $2n+1$-dimensional partially complex
(and by this fact more properly called a $3n+2$-dimensional)
determining system using one free parameter.
\begin{equation}\label{determ_hopf}
\begin{aligned}
0&=f(x^*,x^*,\ldots,x^*,\eta)\\
0&=\Delta(x^*,\eta,\i\omega)v\\
0&=c^\H v-1\\
\end{aligned}
\end{equation}

\paragraph{Periodic solutions}
Periodic solutions are found as solutions $(u(s),\,s\in[0,1];T\in\RR)$
of the following $(n(Ld+1)+1$-dimensional system with no free
parameters.
\begin{align*}\label{determ_psol}
\dot{u}(c_{i,j})=&\
Tf\left(u(c_{i,j}),u\left(\left[c_{i,j}-\frac{\tau_1}{T}\right\vert_{\mod[0,1]}\right),\ldots,
u\left(\left[c_{i,j}-\frac{\tau_m}{T}\right\vert_{\mod[0,1]}\right)\right)\mbox{,}\\
&\quad i=0,\ldots,L-1,\ j=1,\ldots,d \numberthis\\
u(0)=&\ u(1), \\
p(u)=&\ 0.
\end{align*}
Here the notation $t\vert_{\mod[0,1]}$ refers to
$t-\max\{k\in\ZZ:k\leq t\}$, and $p$ represents the integral phase
condition
\begin{equation}\label{integral_phase_cond}
\int_0^1\dot{u}(s)\Delta u(s)\d s=0,
\end{equation}
where $u$ is the current solution and $\Delta u$ its correction.
The collocation points are obtained as 
\[
c_{i,j}=t_i+c_j(t_{i+1}-t_i),\ i=0,\ldots,L-1,\ j=1,\ldots,d,
\]
from the interval points $t_i$, $i=0,\ldots,L-1$ and the collocation
parameters $c_j$, $j=1,\ldots,d$.  The profile $u$ is discretized as a
piecewise polynomial as explained in section \ref{data_structures}.
This representation has a discontinuous derivative at the interval
points. If $c_{i,j}$ coincides with $t_i$ the right derivative is
taken in (\ref{determ_psol}), if it coincides with $t_{i+1}$ the left
derivative is taken. In other words the derivative taken at $c_{i,j}$
is that of $u$ restricted to $[t_i,t_{i+1}]$.

\paragraph{Connecting orbits}
Connecting orbits can be found as solutions of the following
determining system with $s^+-s^-+1$ free parameters, where $s^+$ and
$s^-$ denote the number of unstable eigenvalues of $x^+$ and $x^-$
respectively.
{\allowdisplaybreaks
\begin{align*}
  \dot{u}(c_{i,j})=& Tf(u(c_{i,j}),u(c_{i,j}-\frac{\tau_1}{T}),\ldots,
  u(c_{i,j}-\frac{\tau_m}{T}),\eta)=0\mbox{,\quad}(i=0,\ldots,L-1,\ j=1,\ldots,d) \\
  u(\tilde{c})=&x^{-}+\epsilon
  \sum_{k=1}^{s^{-}}\alpha_{k}v_{k}^{-}e^{\lambda_{k}^{-}T\tilde{c}}, \qquad \tilde{c}<0\\
  0=&f(x^{-},x^{-},\eta)\\
  0=&f(x^{+},x^{+},\eta)\numberthis\label{determ_hcli}\\
  0=&\Delta(x^{-},\lambda_k^-,\eta) v_{k}^{-} \\
  0=&c_k^{H}v_{k}^{-}-1\mbox{,\quad} (k=1,\ldots,s^{-})\\
  0=&\Delta^{H}(x^{+},\lambda_k^+,\eta)w_{k}^{+} \\
  0=&d_k^{H}w_{k}^{+}-1\mbox{,\quad} (k=1,\ldots,s^{+})\\
  0=&{w_{k}^{2}}^{H}(u(1)-x^{+})+\sum_{i=1}^{G}g_i{w_{k}^{+}}^{H}e^{-\lambda_{k}^{+}
    (\theta_i+\tau)}A_{1}(x^{+},\eta)
  \left(u(1+\frac{\theta_i}{T})-x^{+}\right)\mbox{,\quad}(k=1,\ldots,s^+)\\
  u(0)=&x^{-}+\epsilon\sum_{i=1}^{s^{-}}\alpha_{k}v_{k}^{-}\\
  1=&\sum_{i=1}^{s^-}|\alpha_{k}|^{2}\\
  0=&p(u,\eta)
\end{align*}
}Again, all arguments of $u$ are taken modulo $[0,1]$. We choose the
same phase condition as for periodic solutions and similar
normalization of $v_k^-$ and $w+k^+$ as in \eqref{determ_hopf}.

\paragraph{Point method parameters}
The point method parameters (see table \ref{point_method_structures}) 
specify the following options.
\begin{itemize}
\item \blist{newton_max_iterations}: maximum number of Newton
iterations.    
\item \blist{newton_nmon_iterations}: during a first phase of
  \blist{newton_nmon_iterations}+1 Newton iterations the norm of the
  residual is allowed to increase. After these iterations, corrections
  are halted upon residual increase.
\item \blist{halting_accuracy}: corrections are halted when the norm
  of the last computed residual is less than or equal to
  \blist{halting_accuracy} is reached.
\item \blist{minimal_accuracy}: a corrected point is accepted when the
  norm of the last computed residual is less than or equal to
  \blist{minimal_accuracy}.
\item \blist{extra_condition}: this parameter is nonzero when extra
  conditions are provided in a routine \file{sys\_cond.m} which should
  border the determining systems during corrections.  The routine
  accepts the current point as input and produces an array of
  condition residuals and corresponding condition derivatives (as an
  array of point structures) as illustrated below
  (\S\ref{extra_cond}).
\item \blist{print_residual_info}: when nonzero, the Newton iteration
  number and resulting norm of the residual are printed to the screen
  during corrections.
\end{itemize}
For periodic solutions and connecting orbits, the extra mesh
parameters (see table \ref{point_method_structures}) provide the following
information.
\begin{itemize}
\item \blist{phase_condition}: when nonzero the integral phase
  condition \eqref{integral_phase_cond} is used.
\item \blist{collocation_parameters}: this parameter contains user
  given collocation parameters. When empty, Gauss-Legendre collocation
  points are chosen.
\item \blist{adapt_mesh_before_correct}: before correction and if the
  mesh inside the point is nonempty, adapt the mesh every
  \blist{adapt_mesh_before_correct} points.  E.g.: if zero, do not
  adapt; if one, adapt every point; if two adapt the points with odd
  point number.
\item \blist{adapt_mesh_after_correct}: similar to
  \blist{adapt_mesh_before_correct} but adapt mesh after successful
  corrections and correct again.
\end{itemize}

\subsection{Extra conditions}\label{extra_cond}

When correcting a point or computing a branch, it is possible to add
one or more extra conditions and corresponding free parameters to the
determining systems presented earlier. These extra conditions should
be implemented using a function \blist{sys_cond} and setting the
method parameter \blist{extra_condition} to \blist{1} (cf.\ table
\ref{point_method_structures}).  The function \blist{sys_cond} accepts
the current point as input and produces a residual and corresponding
condition derivatives (as a point structure) per extra condition.

As an example, suppose we want to compute a branch of
periodic solutions of system \eqref{example_sys} subject to the following
extra conditions 
\begin{equation}\label{eq:extra_cond}
  \begin{split}
    T&=200, \\
    0&=a_{12}^2+a_{21}^2-1\mbox{,}
  \end{split}
\end{equation}
that is, we wish to continue a branch with fixed period $T=200$ and
parameter dependence $a_{12}^2+a_{21}^2=1$.  The routine shown in
Listing~\ref{sys_cond_demo} implements these conditions by evaluating
and returning each residual for the given point and the derivatives of
the conditions w.r.t.\ all unknowns (that is, w.r.t.\ to all the
components of the point structure).
\begin{lstlisting}[frame=lines,label=sys_cond_demo,caption={Implementation
    extra conditions \eqref{eq:extra_cond} using a routine \blist{sys_cond}.}]
function [resi,condi]=sys_cond(point)
% kappa beta a12 a21 tau1 tau2 tau_s
if point.kind=='psol'
  % fix period at 200:
  resi(1)=point.period-200;
  % derivative of first condition wrt unknowns:
  condi(1)=p_axpy(0,point,[]);
  condi(1).period=1;
  % parameter condition:
  resi(2)=point.parameter(3)^2+point.parameter(4)^2-1;
  % derivative of second condition wrt unknowns:
  condi(2)=p_axpy(0,point,[]);
  condi(2).parameter(3)=2*point.parameter(3);
  condi(2).parameter(4)=2*point.parameter(4);
else
  error('SYS_COND: point is not psol.');
end
end
\end{lstlisting}

\subsection{Continuation}\label{continuation}

During continuation, a branch is extended by a combination of
predictions and corrections.  A new point is predicted based on
previously computed points using secant prediction over an appropriate
steplength. The prediction is then corrected using the determining
systems \eqref{determ_stst}, \eqref{determ_fold}, \eqref{determ_hopf},
\eqref{determ_psol} or \eqref{determ_hcli} bordered with a steplength
condition which requires orthogonality of the correction to the secant
vector.  Hence one extra free parameter is necessary compared to the
numbers mentioned in the previous section.
 
The following continuation and steplength determination strategy is
used.  If the last point was successfully computed, the steplength is
multiplied with a given, constant factor greater than 1.  If
corrections diverged or if the corrected point was rejected because
its accuracy was not acceptable, a new point is predicted, using
linear interpolation, halfway between the last two successfully
computed branch points.  If the correction of this point succeeds, it
is inserted in the point array of the branch (before the previously
last computed point).  If the correction of the interpolated point
fails again, the last successfully computed branch point is rejected
(for fear of branch switch) and the interpolation procedure is
repeated between the (new) last two branch points. Hence, if, after a
failure, the interpolation procedure succeeds, the steplength is
approximately divided by a factor two. Test results indicate that this
procedure is quite effective and proves an efficient alternative to
using only (secant) extrapolation with steplength control.  The reason
for this is mainly that the secant extrapolation direction is not
influenced by halving the steplength but it is by inserting a newly
computed point in between the last two computed points.

\paragraph{Continuation method parameters}
The continuation method parameters (see table \ref{continuation_structure})
have the following meaning.
\begin{itemize}
\item \blist{plot}: if nonzero, plot predictions and corrections
during continuation.
\item \blist{prediction}: this parameter should be 1, indicating
that secant prediction is used (being currently the only
alternative).
\item \blist{steplength_growth_factor}: grow the steplength with
this factor in every step except during interpolation.
\item \blist{plot_progress}: if nonzero, plotting is visible during
  continuation process. If zero, only the final result is drawn.
\item \blist{plot_measure}: if empty use default measures to plot.
  Otherwise \blist{plot_measure} contains two fields, 'x' and 'y',
  which contain measures (see table \ref{measure_structure}) for use
  in plotting during continuation.
\item \blist{halt_before_reject}: If this parameter is nonzero,
continuation is halted whenever (and instead of) rejecting a
previously accepted point based on the above strategy.
\end{itemize}

\subsection{Roots of the characteristic equation}\label{root_char_equa_gio_label}
Roots of the characteristic equation are approximated using a linear
multi-step (LMS-) method applied to \eqref{the_var_equa}.

Consider the linear $k$-step formula
\begin{equation}\label{lms_method}
\sum_{j=0}^k\alpha_j y_{L+j}=h\sum_{j=0}^k\beta_j f_{L+j}.
\end{equation}
Here, $\alpha_0=1$, $h$ is a (fixed) step size and 
$y_j$ presents the numerical approximation of $y(t)$ at the mesh
point $t_j\defeq jh$.
The right hand side
$f_j\defeq f(y_j,\tilde{y}(t_j-\tau_1),\ldots,\tilde{y}(t_j-\tau_m))$ 
is computed using approximations $\tilde{y}(t_j-\tau_1)$ 
obtained from $y_i$ in the past, $i<j$.
In particular, the use of so-called Nordsieck interpolation, leads to
\begin{equation}\label{past_terms}
\tilde{y}(t_j+\epsilon h)=\sum_{l=-r}^s P_l(\epsilon)y_{j+l},\ \epsilon \in [0,1).
\end{equation}
using
\[
P_l(\epsilon)\defeq\prod_{k=-r,\,k\neq l}^s\frac{\epsilon-k}{l-k}.
\]
The resulting method is explicit whenever $\beta_0=0$ and 
$\min{\tau_i}>sh$.
That is, $y_{L+k}$ can then directly be computed from (\ref{lms_method})
by evaluating
\[
y_{L+k}=-\sum_{j=0}^{k-1}\alpha_j y_{L+j}+h\sum_{j=0}^k\beta_j f_{L+j}.
\]
whose right hand side depends only on $y_j$, $j<L+k$.

For the linear variational equation (\ref{the_var_equa})
around a steady state solution $x^*(t)\equiv x^*$
we have
\begin{equation}\label{linear_rhs}
f_j=A_0y_j+\sum_{i=0}^mA_i\tilde{y}(t_j-\tau_i)
\end{equation}
where we have omitted the dependency of $A_i$ on $x^*$.
The stability of the difference scheme (\ref{lms_method}), (\ref{linear_rhs})
can be evaluated by setting $y_j=\mu^{j-L_{\min}}$, $j=L_{\min},\ldots,L+k$ 
where $L_{\min}$ is the 
smallest index used, taking the determinant of (\ref{lms_method})
and computing the roots $\mu$. If the roots of the
polynomial in $\mu$ all have modulus smaller than unity, the trajectories
of the LMS-method converge to zero. 
If roots exist with modulus greater than unity then trajectories exist
which grow unbounded.

Since the LMS-method forms an approximation of
the time integration operator over the time step $h$, so do the 
roots $\mu$ approximate the eigenvalues of $S(h,0)$.
The eigenvalues of $S(h,0)$ are exponential transforms of
the roots $\lambda$ of the characteristic 
equation (\ref{the_char_eq}),
\[
\mu=\exp(\lambda h).
\]
Hence, once $\mu$ is found, $\lambda$ can be extracted using,
\begin{equation}\label{extract_real_part}
\Re(\lambda)=\frac{\ln(|\mu|)}{h}.
\end{equation}
The imaginary part of $\lambda$ is found modulo $\pi/h$, using
\begin{equation}\label{extract_imag_part}
\Im(\lambda)\equiv\frac{\arcsin(\frac{\Im(\mu)}{|\mu|})}
{h}\!\!\!\!\pmod{\frac{\pi}{h}}.
\end{equation}
For small $h$, $0<h\ll 1$, the smallest representation 
in (\ref{extract_imag_part})
is assumed the most accurate one (that is, we let $\arcsin$
map into $[-\pi/2,\pi/2]$).

The parameters $r$ and $s$ (from formula (\ref{past_terms}))
are chosen such that $r\leq s\leq r+2$ (see \cite{Hong96}).
The choice of $h$ is based on the related 
heuristic outlined in \cite{engel01}.

Approximations for the rightmost roots $\lambda$ obtained
from the LMS-method using (\ref{extract_real_part}), 
(\ref{extract_imag_part}) can be corrected
using a Newton process on the system,
\begin{equation}\label{determ_root}
\begin{aligned}
0&=\Delta(\lambda)v \\
0&=c^\T v-1
\end{aligned}
\end{equation}
A starting value for $v$ is the eigenvector of 
$\Delta(\lambda)$ corresponding to its smallest eigenvalue (in modulus).

Note that the collection of successfully corrected roots presents more
accurate yet less robust information than the set of uncorrected
roots. Indeed, attraction domains of roots of equations like
(\ref{determ_root}) can be very small and hence corrections may
diverge, or different roots may be corrected to a single 'exact' root
thereby missing part of the spectrum.  The latter does not occur when
computing the (full) spectrum of a discretization of $S(h,0)$.

Stability information is kept in the structure of table
\ref{stab_structures} (left). The time step used is kept in field
\blist{h}. Approximate roots are kept in field \blist{l0}, corrected
roots in field \blist{l1}.  If unconverged corrected roots are
discarded, field \blist{n1} is empty.  Otherwise, the number of Newton
iterations used is kept for each root in the corresponding position of
\blist{n1}. Here, $-1$ signals that convergence to the required
accuracy was not reached.  

\paragraph{Stability method parameters}
The stability method parameters (see table \ref{meth_stab_struct}
(top)) now have the following meaning.
\begin{itemize}
\item \blist{lms_parameter_alpha}: LMS-method parameters $\alpha_j$ 
ordered from past to present, $j=0,1,\ldots,k$.
\item \blist{lms_parameter_beta}: LMS-method parameters $\beta_j$ 
ordered from past to present, $j=0,1,\ldots,k$.
\item \blist{lms_parameter_rho}: safety radius
  $\rho_{\mathrm{LMS},\epsilon}$ of the LMS-method stability region.
  For a precise definition, see \cite[\S III.3.2]{Enge00}.
\item \blist{interpolation_order}: order of the interpolation in the
past, $r+s=\blist{interpolation_order}$.
\item \blist{minimal_time_step}: minimal time step relative to maximal
delay, $\frac{h}{\tau}\geq\blist{minimal_time_step}$.
\item \blist{maximal_time_step}: maximal time step relative to maximal
delay, $\frac{h}{\tau}\leq\blist{minimal_time_step}$. 
\item \blist{max_number_of_eigenvalues}: maximum number of rightmost
eigenvalues to keep.
\item \blist{minimal_real_part}: choose $h$ such as the discretized
  system approximates eigenvalues with $\Re(\lambda)\geq
  \blist{minimal_real_part}$ well, discard eigenvalues with
  $\Re(\lambda)<\blist{minimal_real_part}$.  If $h$ is smaller than
  its minimal value, it is set to the minimal value and a warning is
  given. If it is larger than its maximal value it is reduced to that
  number without warning.  If minimal and maximal value coincide, $h$
  is set to this value without warning.  If \blist{minimal_real_part}
  is empty, the value $\blist{minimal_real_part}=\frac{1}{\tau}$ is
  used.
\item \blist{max_newton_iterations}: maximum number of Newton
  iterations during the correction process (\ref{determ_root}).
\item \blist{root_accuracy}: required accuracy of the norm of the
  residual of (\ref{determ_root}) during corrections.
\item \blist{remove_unconverged_roots}: if this parameter is zero,
  unconverged roots are discarded (and stability field \blist{n1} is
  empty).
\item \blist{delay_accuracy} (only for state-dependent delays): if the
  value of a state-dependent delay is less than
  \blist{delay_accuracy}, the stability is not computed.
\end{itemize}

\subsection{Floquet multipliers}
Floquet multipliers are computed as eigenvalues of the discretized
time integration operator $S(T,0)$.  The discretization is obtained
using the collocation equations \eqref{determ_psol} without the modulo
operation (and without phase and periodicity condition).  From this
system a discrete, linear map is obtained between the variables
presenting the segment $[-\tau/T,0]$ and those presenting the segment
$[-\tau/T+1,1]$.  If these variables overlap, part of the map is just
a time shift.

Stability information is kept in the structure of table
\ref{stab_structures} (right). Approximations to the Floquet
multipliers are kept in field \blist{mu}.
\paragraph{Stability method parameters}
The stability method parameters (see table \ref{meth_stab_struct}
(bottom)) have the following meaning.
\begin{itemize}
\item \blist{collocation_parameters}: user given collocation
  parameters or empty for Gauss-Legendre collocation points.
\item \blist{max_number_of_eigenvalues}: maximum number of multipliers
  to keep.
\item \blist{minimal_modulus}: discard multipliers with
$|\mu|<\blist{minimal_modulus}$.
\item \blist{delay_accuracy} (only for state-dependent delays): if the
  value of a state-dependent delay is less than
  \blist{delay_accuracy}, the stability is not computed.
\end{itemize}

\section{Concluding comments}\label{limits_sec}

The first aim of {\DDEBIFCODE} is to provide a portable,
user-friendly 
tool for numerical bifurcation analysis 
of steady state solutions and periodic solutions of systems
of delay differential equations of the kinds (\ref{the_dde_type})
and (\ref{the_dde_type2}).
Part of this goal was fulfilled through choosing
the portable, programmer-friendly environment
offered by Matlab.
Robustness with respect to the numerical approximation
is achieved through automatic steplength
selection in approximating the rightmost characteristic roots
and through collocation using piecewise polynomials combined
with adaptive mesh selection.

Although
the package has been successfully tested on a number of realistic examples,
a word of caution may be appropriate. First of all, the package
is essentially a research code (hence we accept no
reliability) in a quite unexplored area of current
research. In our experience up to now, 
new examples did not fail to produce
interesting theoretical questions (e.g., concerning homoclinic
or heteroclinic solutions) many of which
remain unsolved today. 
Unlike for ordinary differential equations, discretization
of the state space is unavoidable during computations on
delay equations. Hence the user of the package is 
strongly advised 
to investigate the effect of discretization using tests on different
meshes and with different method parameters; and, if
possible, to compare with analytical results and/or results obtained
using simulation.

Although there are no 'hard' limits programmed in the package (with
respect to system and/or mesh sizes), the user will notice the rapidly
increasing computation time for increasing system dimension and mesh
sizes.  This is most notable in the stability and periodic solution
computations.  Indeed, eigenvalues are computed from large sparse
matrices without exploiting sparseness and the Newton procedure for
periodic solutions is implemented using direct methods.  Nevertheless
the current version is sufficient to perform bifurcation analysis of
systems with reasonable properties in reasonable execution times.
Furthermore, we hope future versions will include routines which scale
better with the size of the problem.

\subsection{Existing extensions}
\label{sec:extensions}
\begin{itemize}
\item Extension \texttt{debiftool\_extra\_psol} continues the three
  local co\-dim\-ension-one bifurcations of periodic orbits, the fold
  bifurcation, the period doubling and the torus bifurcation for
  constant and state-dependent delays.
\item Extension \texttt{debiftool\_extra\_rotsym} continues relative
  equilibria and relative periodic orbits and their local
  codimension-one bifurcations for constant delays in systems with
  rotational symmetry (that is, there exists a matrix
  $A\in\RR^{n\times n}$ such that $A^T=-A$ and
  $\exp(At)f(x_0,\ldots,x_m)=f(\exp(At)x_0,\ldots,\exp(At)x_m)$).
\item Extension \texttt{debiftool\_extra\_nmfm} computes normal form
  coefficients of Hopf bifurcations, Hopf-Hopf interactions,
  generalized Hopf (Bautin) bifurcations, and zero-Hopf interactions
  (Gavrilov-Guckenheimer bifurcations) for equations with constant
  delays.
\end{itemize}

Other possible future developments include
\begin{itemize}
\item a graphical user interface;
\item incorporation of the numerical core
  routines into a general continuation framework such as \texttt{Coco}
  \cite{DS13} (which would permit the user to grow higher-dimensional
  solution families and wrap other continuation algorithms around the
  core DDE routines),
\item the extension to other types of delay equations such as
  distributed delay and neutral functional differential equations. See
  also Barton \emph{et al} \cite{Barton06} for a demonstration of how
  to extend \DDEBIFCODE{} to neutral functional differential
  equations.
\item determination of more normal-form coefficients to detect
  other co-dimension-two bifurcations.
\end{itemize}
\section*{Acknowledgements}

\DDEBIFCODE{} v.~2.03 is a result of the research project OT/98/16,
funded by the Research Council K.U.Leuven; of the research project
G.0270.00 funded by the Fund for Scientific Research - Flanders
(Belgium) and of the research project IUAP P4/02 funded by the
programme on Interuniversity Poles of Attraction, initiated by the
Belgian State, Prime Minister's Office for Science, Technology and
Culture.  K.~Engelborghs is a Postdoctoral Fellow of the Fund for
Scientific Research - Flanders (Belgium).  J.~Sieber's contribution to
the revision leading to version 3.0 was supported by EPSRC grant
EP/J010820/1.

\bibliographystyle{plain}

\bibliography{manual}

\newpage
\appendix

\section{Jacobians of tutorial examples \blist{neuron} and \blist{sd_demo}}
\label{sec:sys:deri}
The meaning of input arguments to \blist{sys_deri} is explained in
section~\ref{sec:constjac}. The analysis of tutorial example
\blist{neuron} is demonstrated in demo \texttt{neuron} step by step
(see \demobase{}).

\begin{lstlisting}[frame=lines,caption={Jacobians of right-hand side
  \blist{neuron_sys_rhs} in section~\ref{sec:constrhs} for
  \blist{neuron_sys_rhs}.},label=neuron:sys:deri]
function J=neuron_sys_deri(xx,par,nx,np,v)
%% Parameter vector
% par=[\kappa, \beta, a_{12}, a_{21},\tau_1,\tau_2, \tau_s].
J=[];
if length(nx)==1 && isempty(np) && isempty(v)
    %% First-order derivatives wrt to state nx+1
    if nx==0 % derivative wrt x(t)
        J(1,1)=-par(1);
        J(2,2)=-par(1);
    elseif nx==1 % derivative wrt x(t-tau1)
        J(2,1)=par(4)*(1-tanh(xx(1,2))^2);
        J(2,2)=0;
    elseif nx==2 % derivative wrt x(t-tau2)
        J(1,2)=par(3)*(1-tanh(xx(2,3))^2);
        J(2,2)=0;
    elseif nx==3 % derivative wrt x(t-tau_s)
        J(1,1)=par(2)*(1-tanh(xx(1,4))^2);
        J(2,2)=par(2)*(1-tanh(xx(2,4))^2);
    end
elseif isempty(nx) && length(np)==1 && isempty(v)
    %% First order derivatives wrt parameters
    if np==1 % derivative wrt kappa
        J(1,1)=-xx(1,1);
        J(2,1)=-xx(2,1);
    elseif np==2 % derivative wrt beta
        J(1,1)=tanh(xx(1,4));
        J(2,1)=tanh(xx(2,4));
    elseif np==3 % derivative wrt a12
        J(1,1)=tanh(xx(2,3));
        J(2,1)=0;
    elseif np==4 % derivative wrt a21
        J(2,1)=tanh(xx(1,2));
    elseif np==5 || np==6 || np==7 % derivative wrt tau
        J=zeros(2,1);
    end;
elseif length(nx)==1 && length(np)==1 && isempty(v)
    %% Mixed state, parameter derivatives
    if nx==0 % derivative wrt x(t)
        if np==1 % derivative wrt beta
            J(1,1)=-1;
            J(2,2)=-1;
        else
            J=zeros(2);
        end;
    elseif nx==1 % derivative wrt x(t-tau1)
        if np==4 % derivative wrt a21
            J(2,1)=1-tanh(xx(1,2))^2;
            J(2,2)=0;
        else
            J=zeros(2);
        end;
    elseif nx==2 % derivative wrt x(t-tau2)
        if np==3 % derivative wrt a12
            J(1,2)=1-tanh(xx(2,3))^2;
            J(2,2)=0;
        else
            J=zeros(2);
        end;
    elseif nx==3 % derivative wrt x(t-tau_s)
        if np==2 % derivative wrt beta
            J(1,1)=1-tanh(xx(1,4))^2;
            J(2,2)=1-tanh(xx(2,4))^2;
        else
            J=zeros(2);
        end;
    end;
elseif length(nx)==2 && isempty(np) && ~isempty(v)
    %% Second order derivatives wrt state variables
    if nx(1)==0 % first derivative wrt x(t)
        J=zeros(2);
    elseif nx(1)==1 % first derivative wrt x(t-tau1)
        if nx(2)==1
            th=tanh(xx(1,2));
            J(2,1)=-2*par(4)*th*(1-th*th)*v(1);
            J(2,2)=0;
        else
            J=zeros(2);
        end;
    elseif nx(1)==2 % derivative wrt x(t-tau2)
        if nx(2)==2
            th=tanh(xx(2,3));
            J(1,2)=-2*par(3)*th*(1-th*th)*v(2);
            J(2,2)=0;
        else
            J=zeros(2);
        end
    elseif nx(1)==3 % derivative wrt x(t-tau_s)
        if nx(2)==3
            th1=tanh(xx(1,4));
            J(1,1)=-2*par(2)*th1*(1-th1*th1)*v(1);
            th2=tanh(xx(2,4));
            J(2,2)=-2*par(2)*th2*(1-th2*th2)*v(2);
        else
            J=zeros(2);
        end
    end
end

% Raise error if the requested derivative does not exist
if isempty(J)
    error(['SYS_DERI: requested derivative nx=%d, np=%d,',...
        ' size(v)=%d could not be computed!'],nx,np,size(v));
end
end
\end{lstlisting}

The meaning of input arguments to \blist{sys_dtau} is explained in
section~\ref{sec:sddtau}. The analysis of tutorial example
\blist{sd_demo} is demonstrated step by step in demo \texttt{sd\_demo}
(see \demobase{}).
\begin{lstlisting}[frame=lines,caption={Jacobians of the delay function
  \blist{sd_tau} in Listing~\ref{sd_tau} for sd-DDE tutorial example
  \eqref{example_sys2}.},label=sd_dtau]
function dtau=sd_dtau(ind,xx,par,nx,np)
% p1 p2 p3 p4 p5 p6 p7 p8 p9 p10 p11
dtau=[];
if length(nx)==1 && isempty(np)
    %% first order derivatives wrt state variables:
    if nx==0 % derivative wrt x(t)
        if ind==3
            dtau(1:5)=0;
            dtau(2)=par(5)*par(10)*xx(2,2);
        elseif ind==4
            dtau(1)=xx(2,3)/(1+xx(1,1)*xx(2,3))^2;
            dtau(2:5)=0;
        elseif ind==5
            dtau(1:5)=0;
            dtau(4)=1;
        elseif ind==6
            dtau(5)=1;
        else
            dtau(1:5)=0;
        end;
    elseif nx==1 % derivative wrt x(t-tau1)
        if ind==3
            dtau(1:5)=0;
            dtau(2)=par(5)*par(10)*xx(2,1);
        else
            dtau(1:5)=0;
        end;
    elseif nx==2 % derivative wrt x(t-tau2)
        if ind==4
            dtau(1:5)=0;
            dtau(2)=xx(1,1)/(1+xx(1,1)*xx(2,3))^2;
        else
            dtau(1:5)=0;
        end;
    else
        dtau(1:5)=0;
    end;
elseif isempty(nx) && length(np)==1,
    %% First order derivatives wrt parameters
    if ind==1 && np==10
        dtau=1;
    elseif ind==2 && np==11
        dtau=1;
    elseif ind==3 && np==5
        dtau=par(10)*xx(2,1)*xx(2,2);
    elseif ind==3 && np==10
        dtau=par(5)*xx(2,1)*xx(2,2);
    else
        dtau=0;
    end;
elseif length(nx)==2 && isempty(np),
    %% Second order derivatives wrt state variables
    dtau=zeros(5);
    if ind==3
        if (nx(1)==0 && nx(2)==1) || (nx(1)==1 && nx(2)==0)
            dtau(2,2)=par(5)*par(10);
        end
    elseif ind==4
        if nx(1)==0 && nx(2)==0
            dtau(1,1)=-2*xx(2,3)*xx(2,3)/(1+xx(1,1)*xx(2,2))^3;
        elseif nx(1)==0 && nx(2)==2
            dtau(1,2)=(1-xx(1,1)*xx(2,3))/(1+xx(1,1)*xx(2,2))^3;
        elseif nx(1)==2 && nx(2)==0
            dtau(2,1)=(1-xx(1,1)*xx(2,3))/(1+xx(1,1)*xx(2,2))^3;
        elseif nx(1)==2 && nx(2)==2
            dtau(2,2)=-2*xx(1,1)*xx(1,1)/(1+xx(1,1)*xx(2,2))^3;
        end
    end
elseif length(nx)==1 && length(np)==1,
    %% Mixed state parameter derivatives
    dtau(1:5)=0;
    if ind==3
        if nx==0 && np==5
            dtau(2)=par(10)*xx(2,2);
        elseif nx==0 && np==10
            dtau(2)=par(5)*xx(2,2);
        elseif nx==1 && np==5
            dtau(2)=par(10)*xx(2,1);
        elseif nx==1 && np==10
            dtau(2)=par(5)*xx(2,1);
        end
    end
end
%% Otherwise raise error
% Raise error if the requested derivative does not exist
if isempty(dtau)
    error(['SYS_DTAU, delay %d: requested derivative ',...
        'nx=%d, np=%d does not exist!'],ind, nx, np);
end
end

\end{lstlisting} \newpage

\section{Octave compatibility  considerations}
\label{sec:octave}
\paragraph{\blist{nargin} incompatibility} The core DDE-BifTool code
of version v2.0x was likely \texttt{octave} compatible. The changes to \version{},
replacing function names by function handles, broke this compatibility
initially, because, for example, the call \blist{nargin(sys_tau)}
gives an error message in \texttt{octave} (version 3.2.3) if
\blist{sys_tau} is a function handle. To remedy this problem the
additional field \blist{tp_del} is attached to the structure
\blist{funcs} defining the problem. The field \blist{funcs.tp_del} is
set in \blist{set_funcs} (see Section~\ref{sec:funcs:struct} and
Table~\ref{tab:funcs}).

As of version 3.8.1, \texttt{octave} also gives an error message when
one loads function handles as created using \blist{set_funcs} from a
data file. Function handles have to be re-created after a
\blist{clear} or a restart of the session. For an up-to-date list of
known differences in syntax and semantics between matlab and
\texttt{octave} see \url{http://www.gnu.org/software/octave}.
\paragraph{Output} The gradual updating of plots using \blist{drawnow}
slows down for the
\texttt{gnuplot}\footnote{\url{http://www.gnuplot.info/}}-based plot
interface of \texttt{octave} (as of version 3.2.3) as points get added
to the plot. Setting the field \blist{continuation.plot} to
\blist{0.5} (that is, less than $1$ but larger than $0$), prints the
values on the screen instead of updating the plot. The
\texttt{fltk}-based plotting interface of \texttt{octave} does not
appear to experience this slow-down.

Useful options to be set in \texttt{octave}:
\begin{compactitem}
\item \blist{graphics_toolkit('fltk')} sets the graphics toolkit to
  the \texttt{fltk}-based interface (faster plotting);
\item \blist{page_output_immediately(true);} prints out the results of
any \blist{fprintf} or \blist{disp} commands immediately;
\item \blist{page_screen_output(false);} stops paging the terminal
output.
\end{compactitem}

\end{document}